\documentclass[10pt,leqno]{article}

\usepackage{amsmath,amsfonts,amscd,amssymb,theorem}

\long\def\comment#1\endcomment{}

\makeatletter
\begingroup
\gdef\th@dotted{\normalfont\itshape
  \def\@begintheorem##1##2{%
        \item[\hskip\labelsep \theorem@headerfont ##1\ ##2.]}%
\def\@opargbegintheorem##1##2##3{%
   \item[\hskip\labelsep \theorem@headerfont ##1\ ##2\ (##3).]}}
\endgroup
\makeatother

\theoremstyle{dotted}

\newtheorem{theorem}{Theorem}[section]
\newtheorem{lemma}[theorem]{Lemma}

\newtheorem{prop}[theorem]{Proposition}

\makeatletter
\begingroup
\gdef\th@upshape{\normalfont
  \def\@begintheorem##1##2{%
        \item[\hskip\labelsep \theorem@headerfont ##1\ ##2.]}%
\def\@opargbegintheorem##1##2##3{%
   \item[\hskip\labelsep \theorem@headerfont ##1\ ##2\ (##3).]}}
\endgroup
\makeatother

\theoremstyle{upshape}

\newtheorem{defn}[theorem]{Definition}
\newtheorem{remark}[theorem]{Remark}
\newtheorem{example}[theorem]{Example}

\makeatletter
\renewcommand{\subsection}{\@startsection{subsection}{2}{0pt}{-3ex
plus -1ex minus -0.2ex}{-2mm plus -0pt minus
-2pt}{\normalfont\bfseries}} 
\renewcommand{\subsubsection}{\@startsection{subsubsection}{2}{0pt}{-3ex
plus -1ex minus -0.2ex}{-2mm plus -0pt minus
-2pt}{\normalfont\bfseries}} 
\makeatother

\makeatletter
\@addtoreset{equation}{section}
\makeatother

\newcommand{\cntrct}                % contraction with a vector field
{\hspace{2pt}\raisebox{1pt}{\text{$\lrcorner$}}\hspace{2pt}}

\newcommand{\proof}[1][Proof.]{\smallskip\noindent{\em #1}}
\def\endproof{\hfill\ensuremath{\square}\par\medskip}

\renewcommand{\labelenumi}{{\normalfont(\roman{enumi})}}

\def\eqref#1{\thetag{\ref{#1}}}

\let\latexref=\ref
\def\ref#1{{\normalfont{\latexref{#1}}}}

\newcommand{\wt}{\widetilde}
\newcommand{\wh}{\widehat}

\newcommand{\hdot}{{\:\raisebox{3pt}{\text{\circle*{1.5}}}}}
\newcommand{\calo}{{\cal O}}

\newcommand{\kk}{{\sf k}}
\newcommand{\Grp}{\Gamma}
\newcommand{\grp}{\gamma}

\newcommand{\tw}{ {(1)} } 
\newcommand{\fmod}{{\text{\rm -mod}^{\text{{\tt\tiny fg}}}}}

\newcommand{\A}{{\mathcal A}}
\newcommand{\B}{{\mathcal B}}
\newcommand{\CA}{{\mathcal A}}

\newcommand{\E}{{\mathcal E}}
\newcommand{\F}{{\mathcal F}}
\newcommand{\D}{{\mathcal D}}
\newcommand{\G}{{\mathcal G}}

\newcommand{\W}{{\mathcal W}}
\newcommand{\T}{{\mathcal T}}

\newcommand{\X}{{\mathfrak X}}

\newcommand{\bW}{{\sf W}_h}
\newcommand{\bWbar}{{\sf W}}
\newcommand{\bO}{{\sf O}}

\newcommand{\Zet}{{\mathbb Z}}

\newcommand{\R}{{\mathcal R}}

\newcommand{\Lotimes}{\overset{\rm L}{\otimes}}
\newcommand{\imbed}{\hookrightarrow}

\newcommand{\iso}{{\widetilde\longrightarrow}}

\renewcommand{\dim}{\operatorname{\sf dim}}
\newcommand{\codim}{\operatorname{\sf codim}}
\newcommand{\cchar}{\operatorname{\sf char}}
\newcommand{\id}{\operatorname{\sf id}}
\newcommand{\rk}{\operatorname{\sf rk}}
\newcommand{\tr}{\operatorname{\sf tr}}
\newcommand{\Fr}{\operatorname{\sf Fr}}
\newcommand{\FrN}{\operatorname{\it Frob}}
\newcommand{\Spec}{\operatorname{Spec}}
\newcommand{\Proj}{\operatorname{Proj}}
\newcommand{\Coh}{\operatorname{Coh}}

\newcommand{\Hom}{\operatorname{Hom}}

\renewcommand{\Im}{\operatorname{{\sf Im}}}
\newcommand{\End}{\operatorname{End}}
\newcommand{\Ext}{\operatorname{Ext}}
\newcommand{\RHom}{\operatorname{RHom}}
\newcommand{\Aut}{\operatorname{Aut}}

\newcommand{\eend}{\operatorname{{\mathcal E}{\it nd}}}

\newcommand{\Pic}{\operatorname{Pic}}
\newcommand{\Hilb}{\operatorname{Hilb}}

\newcommand{\uHom}{{{\mathcal H}{\it om}}}

\newcommand{\uuHom}{{\sf Hom}}
\newcommand{\bbS}{{\mathbb S}}

\newcommand{\Db}{D^b}

\newcommand{\Gm}{{\mathbb{G}_m}}

\newcommand{\g}{{\mathfrak g}}
\newcommand{\ssp}{{{\mathfrak s}{\mathfrak p}}}
\newcommand{\m}{{\mathfrak m}}

\newcommand{\K}{{\sf K}}

\newcommand{\sllash}{/\!\!\,/}

\title{McKay equivalence for symplectic resolutions of singularities}

\author{Roman Bezrukavnikov\thanks{Partially supported by NSF grant
DMS0071967.} and Dmitry Kaledin\thanks{Partially supported by CRDF
grant RM1-2354-MO02.}}

\begin{document}

\maketitle

\tableofcontents

\section{Introduction.}

Let $\K$ be an algebraically closed field of characteristic $0$, let
$V$ be a finite-dimensional $\K$-vector space equipped with a
non-degenerate skew-symmetric form $\omega\in \Lambda^2(V^*)$, and
let $\Grp\subset Sp(V)$ be a finite subgroup.  Suppose that we are
given a resolution of singularities of the quotient variety
$\pi:X\to V/\Grp$ such that the symplectic form on the smooth part
of $V/\Grp$ extends to a non-degenerate closed $2$-form $\Omega
\in H^0(\Omega^2_X)$. The aim of this paper is to prove

\begin{theorem}\label{main}
  There exists an equivalence of $\calo_V^\Grp$-linear
  triangulated categories $D^b(\Coh(X))\cong D^b(\Coh^\Grp(V))$.
\end{theorem}

A conjecture of this type was first made by M. Reid \cite{reid}; a
more general statement was conjectured by A. Bondal and D. Orlov,
\cite[\S 5]{BO}.

When $\dim(V)=2$ such an equivalence is well-known, \cite{KV},
\cite{SV}; in fact our argument relies on these results. Recently a
similar statement was established by T. Bridgeland, A. King and M.
Reid \cite{BKR} for crepant resolutions of Gorenstein quotients of
vector spaces of dimension $3$. The approach of \cite{BKR} is very
elegant, but unfortunately, it works only in $\dim 3$. Our paper is
an attempt to generalize at least the result, if not the approach to
higher dimensions. We can only treat symplectic quotients, not
arbitrary Gorenstein ones.  Notice though that our additional
assumption on the resolution is not restrictive --- every crepant
resolution $X$ of a symplectic quotient singularity in fact carries
a non-degenerate symplectic form (see e.g.  \cite{Ka}).
 
The proof uses reduction to positive characteristic, and
quantization of the symplectic variety $X_{\kk}$ over a field $\kk$
of characteristic $p>0$. Our method is suggested by the results of
\cite{MR}.

The key ingredient of the proof is a quantization of $X_{\kk}$ whose
global sections coincide with the standard quantization of
$H^0(\calo_X) = H^0(V,\calo_V)^\Grp$. By this we mean a
deformation of the structure sheaf $\calo_X$ to a sheaf of
non-commutative $\kk[[h]]$-algebras $\calo_h(X)$ such that the
algebra of global sections $H^0(X,\calo_h)$ is identified the
subalgebra $\W^\Grp \subset \W$ of $\Grp$-invariant vectors in
the (completed) Weyl algebra $\W$ of the vector space $V$.

It turns out that over the generic point of $\Spec(\kk[[h]])$, the
quantized algebra $\calo_h$ is an Azumaya algebra over $X^\tw$, the
Frobenius twist of $X$; it also yields an Azumaya algebra on
$X^\tw_\kk$. The category of modules over the latter is the category
of coherent sheaves on some gerb over $X^\tw$.

One then argues that the above Azumaya algebra on $X^\tw$ is {\em
derived affine}, i.e. the derived functor of global sections
provides an equivalence between the derived category of sheaves of
modules, and the derived category of modules over its global
sections; this algebra of global sections is identified with the
algebra $\bWbar^\Grp$, where $\bWbar$ is the reduction of the Weyl
algebra at $h=1$.

Furthermore, for large $p$ we have a Morita equivalence between
$\bWbar^\Grp$ and $\bWbar\#\Grp$, the smash-product of $\bWbar$ and
$\Grp$. Thus we get an equivalence between $ D^b(\W\#\Grp\fmod)$ and
the derived category of modules over the Azumaya algebra on $X^\tw$.
The algebra $\bWbar$ is an Azumaya algebra over $V^\tw$; thus,
roughly speaking, the latter equivalence differs from the desired
one by a twist with a certain gerb.  We then use the norm map on
Brauer groups to pass from sheaves over a gerb to coherent sheaves
on the underlying variety.

Then the equivalence over $\kk$ of large positive characteristic is
constructed; by a standard procedure we derive the desired statement
over a field of characteristic zero.

\medskip

Theorem~\ref{main} implies, more or less directly, that any crepant
resolution $X$ of the quotient $V/\Grp$ is the moduli space of
$\Grp$-equivariant Artinian sheaves on $V$ satisfying some stability
conditions (what is known nowadays as {\em
$G$-constellations}). This is an important issue which deserves
further study; however, spelling out precisely the stability
conditions is a non-trivial problem, and we prefer to postpone it to
a future paper. Among other things, in the case when $X =
\Hilb^n({\mathbb A}^2)$ is the Hilbert scheme of $n$ points on the
affine plane, this should give a computation-free proof of the
so-called $n!$-Conjecture established recently by M. Haiman
\cite{Haim}. We also note that once one knows that $X$ is a moduli
space of $G$-constellations, one can prove Theorem~\ref{main} by the
method of \cite{BKR}; this is pretty useless, however, because to
obtain a modular interpretation of an arbitrary smooth crepant
resolution $X$ one already needs Theorem~\ref{main}. We would like
to mention that \cite[Corollary 1.3]{BKR} is completely misleading
in this respect, since it in fact assumes that the resolution is a
$G$-Hilbert scheme (although this is not mentioned in the statement,
and only appears as a sort of carry-over assumption from the
previous page).

\medskip

It is easy to see that the category of $\Grp$-equivariant coherent
sheaves on $V$ is equivalent to the category of finitely-generated
modules over a certain non-commutative algebra. Recently B. Keller
\cite{kel} has given a definition of the so-called {\em Hochschild
homology} groups of an abelian category. The Hochschild homology of
the category of finitely-generated modules over an algebra coincides
with the Hochschild homology of the algebra; in particular, the
homology group $HH_k$ is trivial for $k < 0$. On the other hand, the
Hochschild homology group $HH_k$ of the category of coherent sheaves
on a smooth algebraic variety $X$ over a characteristic $0$ field is
isomorphic to the sum of $H^{p-k}(X,\Omega^p_X)$. But Keller's
definition is invariant with respect to functors which induce an
equivalence of derived categories. Therefore Theorem~\ref{main}
implies that $H^p(X,\Omega^q_X) = 0$ for $p > q$. In the symplectic
case, this reduces to the more familiar $H^p(X,\Omega^q_X) = 0$,
$p+q > \dim X$, which has been proved in \cite{K3} (in fact, we use
this vanishing to prove Theorem~\ref{main}). However, in general,
the required vanishing is very strange and quite strong. This raises
some doubts as to whether Theorem~\ref{main} holds for crepant
resolution of general Calabi-Yau quotient singularities $V/\Grp$.

\subsection{Notations.}\label{nota}
The pair $\langle V,\omega\rangle$ is a symplectic vector space,
$\Grp\subset Sp(V)$ is a finite subgroup, $\pi:X\to V/\Grp$ is a
fixed resolution of singularities such that $\omega$ extends to a
symplectic form on the whole $X$, $\eta: V\to V/\Grp$ is the
projection map. The group algebra of $\Grp$ with coefficients in a
ring $R$ is denoted by $R[\Grp]$. For an arbitrary $\kk$-algebra $A$
equipped with a $\Grp$-action, the {\em smash-product algebra} $A \#
\Grp$ coincides with $A [\Grp]$ as an abelian group, and the
multiplication is defined by
$$
(a_1 \cdot \gamma_1)(a_2 \cdot \gamma_2) = a_1a_2^{\gamma_1} \cdot
\gamma_1\gamma_2, \qquad a_1,a_2 \in A, \gamma_1,\gamma_2 \in
\Grp.
$$
The  Weyl algebra $\bW$ is defined by
$$
\bW=\kk[h]\langle V^*\rangle/\left(
  xy-yx=h\omega^{-1}(x,y)\right);
$$
the formal Weyl algebra $\W$ is the $h$-adic completion of $\bW$.
We define open $\Grp$-invariant subschemes $V_i\subset V$ by
$$
V_i=\{v\in V\ |\ \dim V^{Stab_\Grp(v)}\leq 2i+\dim V^\Grp\},
$$
and we set $X_i=\pi^{-1}(V_i/\Grp)$ (these will actually be used for
$i=0,1$ only). We write $\calo^p=\{f^p\ |\ f\in \calo\}$ for a
commutative algebra $\calo$ of characteristic $p$.  For a scheme $X$
over a field of characteristic $p$, $X^\tw$ is the Frobenius twist
of $X$, and $\Fr:X\to X^\tw$ is the Frobenius morphism. In our
applications $X$ will be reduced; in this case $X^\tw =\langle X,
\calo_X^p \rangle$. Moreover, the base field will be perfect, so
that the twist $X^\tw$ is isomorphic to $X$ as an abstract scheme;
however, it is often convenient to distinguish between the two
notation-wise.
 
If $\A$ is an Azumaya algebra over a scheme $X$, we write
$\Coh(X,\A)$ for the category of coherent sheaves of $\A$-modules.

\subsection{Acknowledgements.} We would like to thank  D.~Arinkin,
M.~Finkelberg, V.~Ginzburg, A.~Kuznetsov and V.~Vologodsky for
friendly help with various mathematical questions (in particular,
Remark~\ref{rem_pro_BKR} was pointed out by Kuznetsov).  The second
author is grateful to Professor M.~Lehn for fruitful discussions.
Finally, we would like to acknowledge the intellectual debt to
I.~Mirkovic and D.~Rumynin. This work, at least partly, stems from
an attempt to isolate the relevant geometric context for their ideas
on geometry and representation theory in positive characteristic.

Several years ago the first author gave a talk at the Shafarevich
seminar at Steklov Math Institute in Moscow, where he described the
results of Mirkovic and Rumynin and suggested, somewhat sceptically,
some possible generalizations; the present paper is in large part
the realization of the program sketched at that time. A.N. Tyurin
was in the audience. He immediately believed in our vague hopes and
enthusiastically supported the proposed research. This paper
confirms that he was right. We were hoping to tell Andrei
Nikolaevich about our results; sadly, this was not to be.

\section{Almost exceptional objects.}

The goal of this section is to show that Theorem~\ref{main} in large
part follows from completely general homological arguments.

\begin{defn}
A nonzero object $M$ of an abelian category is {\it almost
exceptional} if $\Ext^i(M,M)=0 $ for $i > 0$, and the algebra
$\End(M)$ has finite homological dimension.
\end{defn}

\begin{prop}\label{CYProp} 
Let $Y$ be a smooth irreducible variety over a field $k$ with
trivial canonical class; assume that there exists a proper morphism
from $Y$ to an affine $k$-variety $S$. Let $\A$ be an Azumaya
algebra over $Y$. If $\E\in \Coh(X,\A)$ is an almost exceptional
object, then the functor $\F\mapsto\RHom^\hdot(\E,\F)$ from the
derived category $D^b(\Coh(X,\A))$ to the derived category
$D^b(\End(\E)^{op}\fmod)$ of finitely generated right
$\End(\E)$-modules is an equivalence.
\end{prop}

Proposition~\ref{CYProp} applied to $\A=\calo$ immediately shows
that to prove Theorem~\ref{main}, it suffices to prove the
following.

\begin{theorem}\label{evb}
There exists a vector bundle $\E$ on $X$ such that
\begin{enumerate}
\item the $\calo_V^\Grp$-algebra $\End(\E)$ is isomorphic to the
  smash-product $\calo_V \# \Grp$,
\item we have $\Ext^i(\E,\E)=0$ for $i > 0$.
\end{enumerate}
\end{theorem}

This is what we will do in the rest of the paper, after proving
Proposition~\ref{pf.CY} in Subsection~\ref{pf.CY} (to eliminate
confusion between right and left modules, we note right away that
since $V \cong V^*$ as a $\Grp$-module, we have
$(\calo_V\#\Grp)^{op}\cong\calo_V\#\Grp$). Additionally, we will
show in Subsection~\ref{rg} that almost exceptional objects are
remarkably rigid, so that essentially it suffices to construct such
a bundle $\E$ after reduction to positive characteristic.

\begin{remark}
One can show that conditions on the vector bundle $\E$ in Theorem
\ref{evb} are equivalent to the following conditions.
\begin{enumerate}
\renewcommand{\labelenumi}{{\normalfont(\alph{enumi})}}
\item $\E_{X_0}\cong \pi^*\eta_*(\calo)|_{X_0}$

\item $\Ext^i(\E,\E)=0$ for $i>0$, and $\End(\E)\iso
\End(\E|_{X_0})$.
\end{enumerate}
\end{remark}

\subsection{Calabi-Yau categories.}\label{CY}

We recall some generalities from homological algebra. We refer to
\cite{BK} for details.

Let $k$ be a field. A $k$-linear triangulated category is said to be
{\em of finite type} if for any two objects $X,Y\in D$ the space
$\Ext^\hdot(X,Y)= \bigoplus_{i\in \Zet} \Hom(X,Y[i])$ is
finite-dimensional. For such a category Bondal and Kapranov
\cite{BK} defined a {\it Serre functor} to be a pair consisting of a
(covariant) auto-equivalence $S:D\to D$, and an isomorphism
$\Hom(X,Y)\cong \Hom(Y,S(X))^*$ fixed for all $X,Y\in D$, subject to
certain compatibilities. If a Serre functor exists, it is unique up
to a unique isomorphism. For example, if $D$ is the bounded derived
category of coherent sheaves on a smooth projective variety $X$ over
$k$, then $F\mapsto F\otimes K_X [\dim X]$ is a Serre functor for
$D$ (where $K_X$ is the canonical line bundle).

We will need a slight generalization of this notion. Let $\calo$ be
a finitely generated commutative algebra over a field, and let $D$
be an $\calo$-linear triangulated category; moreover, assume given a
functor $\uuHom: D^{op}\times D\to D^b(\calo\fmod)$ and a functorial
isomorphism $\Hom(X,Y)\cong H^0(\uuHom(X,Y))$; we will call such a
set of data {\em a strong $\calo$-category}.

The triangulated category $D^b(\calo\fmod)$ is equipped with a
canonical anti-auto-equivalence, namely, the Grothendieck-Serre
duality \cite{Ha}. We let $\bbS$ denote this functor; thus $\bbS:
\F\mapsto R\uHom(\F,\D)$, where $\D$ is the Grothendieck dualizing
sheaf, and $\uHom$ stands for the internal $Hom$.

\begin{defn}
A {\it Serre functor with respect to $\calo$} (or just an {\em
$\calo$-Serre functor}) is a pair consisting of a (covariant)
auto-equivalence $S:D\to D$, and a natural (functorial) isomorphism
$\uuHom(X,Y)\cong \bbS(\uuHom(Y,S(X)))$ satisfying the
compatibilities in \cite{BK}.
\end{defn}

For example, if $X$ is a smooth variety over $k$ equipped with a
projective morphism $\pi:X\to \Spec(R)$, then $D^b(\Coh(X))$ with
$\uuHom(\F,\G) =R\pi_* \uHom(\F,\G)$ is a strong $R$-category. The
functor $\F\mapsto \F\otimes K_X[\dim X]$ is naturally a Serre
functor with respect to $R$; this is true because Grothendieck-Serre
duality commutes with proper direct images, and we have
$$
\bbS (\uHom(\F,\G))\cong \uHom (\G, \F\otimes K_X[\dim X]).
$$
The following generalization of this fact is straightforward.

\begin{lemma}\label{AzuSerre}
  Let $X$ be a smooth variety over $k$ equipped with a projective
  morphism $\pi:X\to \Spec(R)$ as above; and let $\CA$ be an
  Azumaya algebra on $X$. Then $D^b(\Coh(X,\A))$ is
  naturally a strong $R$-category. The functor $\F\mapsto
  \F\otimes K_X[\dim X]$ is naturally a Serre functor with respect
  to $R$.\endproof
\end{lemma}

A strong $\calo$-category will be called {\it Calabi-Yau} if it
admits a Serre functor with respect to $\calo$ which is isomorphic
to the shift functor $X\mapsto X[n]$ for some $n\in \Zet$.

Application of the above notions to our situation is based on the
following 

\begin{lemma}\label{abstract} Let $F:C\to D$ be a triangulated functor
  between non-zero triangulated categories. Assume that
\begin{enumerate}
\item $F$ has a left adjoint functor $F'$ and the adjointness
  morphism $id\to F\circ F'$ is an isomorphism.
  
\item $C$ is indecomposable, i.e. it can not be written as
  $C=C_1\oplus C_2$ for nonzero triangulated categories $C_1,C_2$.
  
\item $C$ admits a structure of a strong $\calo$-category for some
  commutative algebra $\calo$ of finite type over a field, so that
  $C$ is Calabi-Yau with respect to $\calo$.
\end{enumerate}
Then $F$ is an equivalence.
\end{lemma}

\proof{} Condition \thetag{i} implies that $F'$ is a full embedding,
and its (essential) image (which we denote by $I$) is a right
admissible subcategory in $C$. Recall that this means that every
object $\F$ of $C$ fits into an exact triangle $\F_1\to \F\to
\F_2\to \F_1[1]$ where $\F_1\in I$, and $\F_2\in I^\perp$; here
$$
I^\perp=\left\{\G\in C\ \mid\ Hom(\F,\G)=0\quad \forall \F\in
  I\right\}
$$
is the right orthogonal of $C$ (to get such a triangle for a given
$\F$ set $\F_1=F'F(\F)$, and let the arrow $\F_1\to \F$ be the
adjunction morphism; then complete it to an exact triangle). In this
situation one says that $C$ admits a semi-orthogonal
decomposition. It follows from the definitions that if $S$ is a
Serre functor for $C$ (relative to some commutative algebra $\calo$
such that $C$ is equipped with a strong $\calo$-linear structure)
then $S^{-1}$ sends the right orthogonal to a full subcategory into
the left orthogonal to the same subcategory. In particular, if $C$
is Calabi-Yau relative to some algebra $\calo$, then the left
orthogonal to any triangulated subcategory coincides with the right
orthogonal. Thus the above semi-orthogonal decomposition is actually
orthogonal, i.e. we have $C=I\oplus I^\perp$; since $D\ne 0$
condition (ii) implies that $I^\perp=0$, so $I=C$, and $F$, $F'$ are
equivalences. \endproof

\begin{remark}\label{rem_pro_BKR}
A similar (more general) statement is proved and used in  \cite{BKR}.
\end{remark}

Another simple auxiliary fact is

\begin{lemma}\label{indecompo}
  Let $X$ be a connected quasiprojective variety over a field, and
  let $\CA$ be an Azumaya algebra on $X$. Then the category
  $\Db(\Coh(X,\CA))$ is indecomposable.
\end{lemma}

\proof{} Assume that $\Db(\Coh(X,\A))=D_1\oplus D_2$ is a
decomposition invariant under the shift functor. Let $P$ be an
indecomposable summand of the free $\CA$-module. Let $L$ be an ample
line bundle on $X$ such that $H^0(L\otimes \uHom_\calo (P,P))\ne 0$.
For any $n\in \Zet$ the $\CA$-module $P\otimes L^{\otimes n}$ is
indecomposable, hence belongs either to $D_1$ or to $D_2$.
Moreover, all these modules belong to the same summand, because
$Hom(P\otimes L^{\otimes n}, P\otimes L^{\otimes m})\ne 0$ for
$n\leq m$.  If $\F$ is an object of the other summand, we have
$Ext^\hdot (P\otimes L^{\otimes n}, \F)=0$ for all $n$, which
implies $\F=0$. \endproof

\subsection{Proof of Proposition~\ref{CYProp}.}\label{pf.CY}
We claim that the functor 
$$
F:\F \mapsto \RHom^\hdot(\E,\F)
$$
satisfies the conditions of Lemma~\ref{abstract}. Indeed, the left
adjoint functor is given by $F':M\mapsto M\Lotimes_{\End(\E)^{op}}
\F$; it sends the bounded derived category in the bounded one
because $\End(\E)$ has finite homological dimension.  Vanishing of
$Ext^i(\E,\E)$ implies that $\id\cong F'\circ F$. Thus \thetag{i}
holds. Condition \thetag{ii} is provided by Lemma~\ref{indecompo},
and \thetag{iii} holds by Lemma~\ref{AzuSerre}. \endproof

\subsection{Rigidity.}\label{rg}
Let $R$ be a regular Noetherian ring with a maximal ideal $\m
\subset R$ and the residue field $k = R/\m$. Assume given an
algebraic variety $X_R$ flat and smooth over $R$, and let $X = X_R
\otimes_R k$ be the fiber of $X$ over $\Spec k \in Spec R$.

\begin{lemma}\label{rig}
Every almost exceptional vector bundle $\E$ on $X$ extends uniquely
to an almost exceptional vector bundle $\wh{\E}$ on the formal
completion $\X$ of $X_R$ in $X \subset X_R$.
\end{lemma}

\proof{} Filter $R$ by the powers of $\m$, and construct the
extension step-by-step, by extending to $X_R \otimes_R (R/\m^k)$ for
all $k$ in turn. By the standard deformation theory, at each step
the obstructions to extending $\E$ lie in $\Ext^2(\E,\E)$, and
different extensions are parametrized by a torsor over
$\Ext^1(\E,\E)$. Since $\E$ is almost exceptional, both groups
vanish. It remains to show that the extended bundle $\wh{\E}$ is
almost exceptional. Indeed, since $\Ext^i(\E,\E) = 0$ for all $i
\geq 1$, we also have $\Ext^i(\wh{\E},\wh{\E}) = 0$ for $i \geq 1$,
$\End(\wh{\E})$ is flat over $R$, and the natural map
$\End(\wh{\E})/\m \to \End(\E)$ is an isomorphism. To prove that the
algebra $\End(\wh{\E})$ has finite homological dimension, one
computes the $\Ext^\hdot$-groups by using the spectral sequence
associated to the $\m$-adic filtration.
\endproof

Our varieties will usually be non-compact, so that passing from $\X$
to $X_R$ is not automatic. To obtain global information, we will use
$\Gm$-actions. Assume that $X_R = \Proj\B^\hdot$ is projective over
an affine variety $\Spec \B^0$ over $R$.

\begin{defn}
We will say that an action of the group $\Gm$ on $X_R$ is {\em
positive-weight} if all weights of corresponding $\Gm$-action on
$\B^\hdot$ are non-negative.
\end{defn}

\begin{lemma}\label{glob}
In the assumptions of Lemma~\ref{rig}, assume in addition that $R$
is complete with respect to $\m$-adic filtration, $X_R$ is equipped
with a positive-weight $\Gm$-actions and that $\E$ is
$\Gm$-equivariant. Then $\E$ extends uniquely to a $\Gm$-equivariant
vector bundle $\E_R$ on $X_R$.
\end{lemma}

\proof[Sketch of a proof.] Every extension of the multiplicative
group by a unipotent group is split; it follows that the extension
of $\E$ to $\X$ provided by Lemma~\ref{rig} can be equipped with a
$\Gm$-equivariant structure. Under the positive-weight assumption,
the completion functor is an equivalence of categories between
$\Gm$-equivariant coherent sheaves on $X_R$ and on $\X$. To
construct an inverse equivalence, one uses the Serre Theorem to
interpret coherent sheaves as graded $\B^\hdot$-modules, and
replaces a module with the the submodule of $\Gm$-finite sections.
\endproof

In fact, it is not necessary to require that $\E$ is
$\Gm$-equivariant -- every almost exceptional vector bundle is
automatically $\Gm$-equivariant with respect to any positive-weight
$\Gm$-action. We do not prove this, because the proof is slightly
technical; in our applications, $\Gm$-equivariant structure follows
from a direct geometric argument (see Proposition~\ref{lifts}).

\section{Quantizations.}

In this section we spell out some of the generalities on
quantization of algebraic varieties in positive characteristic. We
do not strive for generality, and only do the work needed for the
proof of Theorem~\ref{main}.

Till the end of the section all objects are assumed to be defined
over a fixed field $\kk$.

\subsection{Standard definitions.}

For a Poisson $\kk$-algebra $\langle\calo,\{-,-\}\rangle$, its
(formal) quantization is defined in the usual way; thus a
quantization is an associative flat $\kk[[h]]$-algebra $\calo_h$,
complete separated in the topology generated by $h^i\calo_h$, and
equipped with an isomorphism $\calo_h/h\calo_h\cong \calo$ such that
the commutator in $\calo_h$ equals $h\{-,-\}\mod h^2\calo_h$.

If $X$ is Poisson scheme, a quantization of $X$ is a sheaf of
$\kk[h]$-algebras $\calo_h$ on $X$ equipped with an isomorphism
$\calo_h/h\calo_h\cong \calo_X$ such that the algebra of sections
$\calo_h(U)$ for an affine open $U$ is a quantization of $\calo(U)$.
One checks that quantizazion of an algebra defines a quantization of
its localization (see e.g. \cite[\S 2.1]{Kapranov}), so that a
quantization of a Poisson algebra $A$ gives a quantization of the
Poisson scheme $X=\Spec A$.

Below we will mostly be concerned with examples when $X$ is a
symplectic variety.

\begin{example}\label{Diff op}
  For a smooth affine variety $M$ over $\kk$, the algebra $D_h(M)$
  of {\it asymptotic differential operators} on $M$ is the
  $h$-completion of the algebra generated by $\calo_M$ and $Vect(M)$
  subject to the usual relations $f_1\cdot f_2=f_1f_2$, $f\cdot
  \xi=f\xi$, $\xi\cdot f-f\cdot\xi=h \xi(f)$,
  $\xi_1\xi_2-\xi_2\xi_1=h[\xi_1,\xi_2]$. As usual one checks easily
  that $D_h(M)$ is a quantization of the sympectic variety $T^*M$.
  
  Gluing the above construction, for any smooth variety $M$ over $\kk$
  one obtains a sheaf $D_h(M)$ which is a quantization of $T^*M$.

Notice that when $\cchar \kk$ is positive, this construction is
related to the so-called ``crystalline'' differential operators in
the terminology of Mirkovic and Rumynin (PD differential operators
in the terminology of \cite{Ogus}) rather than to the more widely
known Grothendieck differential operators (the latter contain
divided powers of a vector field, while the former does not).
\end{example}

We will also need a graded version of quantizations. Assume that a
Poisson scheme $X$ is equipped with a $\Gm$-action such that the
Poisson bracket has weight $-2$ (in other words, $\deg \{f,g\} =
\deg f + \deg g - 2$ for any two homogeneous local functions $f$,
$g$ on $X$). We will say that that a quantization $\calo_h$ of the
scheme $X$ is {\em graded} if it is equipped with a $\Gm$-action
such that $h$ has degree $2$ and the isomorphism $\calo_h/h \cong
\calo_X$ is $\Gm$-equivariant. E.g. the standard quantization $\W$
of a symplectic vector space is graded (see Example~\ref{Weyl
algebra}).

\subsection{Quantizations as deformations.}

By its very definition, quantizations can be studied by deformation
theory. Lately it has become fashionable to reduce all questions of
deformation theory to solving the Maurer-Cartan equation in some
differential graded Lie algebra. In characteristic $0$, one can
conceivably apply this method to quantizations and maybe even obtain
something useful; however, this is irrelevant for us, since we are
mostly interested in the positive characteristic case, where the
differential graded Lie algebra formalism makes no sense. Thus we
have to deal with deformations step-by-step, in the traditional
standard way. This does not give much, but it is possible to prove
one extension result which we will need.

Let $X$ be a smooth variety. We recall that by the standard
deformation theory, deformations of the structure sheaf $\calo_X$ in
the class of sheaves of associative algebras are controlled by the
Hochschild cohomology groups $HH^i(X)$, $i=1,2,3$ --- that is, by
the groups
$$
\Ext^i_{X \times X}(\calo_\Delta,\calo_\Delta), \qquad i =2,3,
$$
where $\Delta \subset X \times X$ is the diagonal. The group
$HH^2(X)$ contains the deformation classes, and the group $HH^3(X)$
contains the obstructions.

\begin{lemma}\label{def}
\begin{enumerate}
\item Let $U \subset X$ be an open subset, and assume that the
restriction map
$$
HH^i(X) \to HH^i(U)
$$
is bijective for $i=1,2$ and injective for $i=3$.Then every
quantization of $U$ extends uniquely to a quantization of
$X$.
\item Assume that $X$ is equipped with a $\Gm$-action which
preserves $U \subset X$, so that the Poisson bracket $\{-,-\}$ in
$\calo_X$ has degree $-2$, and let $HH^\hdot_-(X) \subset
HH^\hdot(X)$, $HH^\hdot_-(U) \subset HH^\hdot(U)$ be the subspaces
of vectors of negative weight with respect to $\Gm$. Assume that the
restriction map $HH^i_-(X) \to HH^i_-(U)$ is bijective for $i=1,2$
and injective for $i=3$. Then every graded quantization of $U$
extends uniquely to a graded quantization of $X$.
\end{enumerate}
\end{lemma}

\proof[Sketch of a proof.] By an order-$n$ deformation $\calo_n$ of
the structure sheaf $\calo_X$ we will understand a sheaf of flat
$\kk[h]/h^{n+1}$-algebras on $X$ equipped with an algebra
isomorphism $\calo_n/h \cong \calo_X$. Fix a deformation
$\calo_{n-1}$ of order $n-1$. Then all order-$n$ deformations
$\calo_n$ of $\calo_X$ equipped with an isomorphism $\calo_n/h^n
\cong \calo_{n-1}$ form a gerb bound by $HH^{\leq 3}(X)$ (this means
that for some complex $C^\hdot$ with cohomology groups $HH^i(X)$,
$i=1,2,3$, and for some fixed element $c \in C^3$, all these
deformations are parametrized by elements $b \in C^2$ with $d(b)=c$,
and isomorphic deformations correspond to homologous elements $b_1$,
$b_2$). Analogously, all such deformations of $\calo_U$ form a gerb
bound by $HH^{\leq 3}(U)$. To prove \thetag{i}, note that the
conditions of the Lemma insure that the restriction induces an
equivalence of the corresponding gerbs, and apply induction on
$n$. To prove \thetag{ii}, note that in the $\Gm$-equivariant
setting, since we have $\deg h =2$, the relevant gerbs are bound by
$HH^{\leq 3}(X)_{-2n} \subset HH^{\leq 3}(X)$, that is, the subspace
of vectors of weight $-2n$. The conditions of the Lemma insure that
the gerbs are equivalent at each step $n$. \endproof

\subsection{Frobenius center.}
From now till the end of the section we assume that $\kk$ is a
perfect field $\kk$ of characteristic $p>0$. The Poisson scheme $X$
is assume to be reduced.

The new feature of the positive characteristic theory is the
presence of a huge center in the quantized algebra $\calo_h$. In the
situation of Example~\ref{Diff op} it is closely related to the
notion of {\it $p$-curvature } of a flat connection.

We start with an elementary observation that for a Poisson
$\kk$-algebra $\calo$ and $f\in \calo$, the element $f^p\in \calo$
lies in the center of the Lie algebra $\calo$; thus the Poisson
bracket is $\calo^p$-linear, so that the Poisson bracket turns
$\Fr_*(\calo)$ into a coherent sheaf of Lie algebras on the scheme
$X^\tw=(X,\calo_X^p)$.

\begin{defn} 
A quantization $\calo_h$ of a Poisson scheme $X$ is called {\it
Frobenius-constant} if the embedding $\calo_X^p\imbed \calo$ lifts
to a morphism of sheaves of algebras from $\calo_X^p$ to the center
of $\calo_h$
\end{defn}

Frobenius-constant quantizations are more geometric than arbitrary
ones. Namely, let $\calo_h$ be a Frobenius-constant quantization of
a Poisson scheme $X$. Then $\calo_h$ is by definition a sheaf of
$\calo^p$-algebras, thus $\calo_h$ defines a quasi-coherent sheaf on
the formal scheme $\widehat{X}^\tw$, the completion of $X^\tw\times
\Spec\kk[h]$ at the central fiber $X^\tw\times \{0\}$. It is easy to
see that this sheaf is locally free of rank $p^{\dim(X)}$. Its
restriction to the special fiber of $\widehat X$ is identified with
$Fr_*(\calo)$. As in Lemma~\ref{glob}, in the graded setting we have
even more.

\begin{lemma}\label{gra}
Let $X = \Proj \B^\hdot$ be a projective variety over $\Spec \B^0$,
smooth over $\kk$ and equipped with a $\Gm$-action and a Poisson
bracket of degree $-2$.  Assume that the $\Gm$-action on $X$ is
positive-weight. Then every graded Frobenius-constant quantization
$\calo_h$ of $X$ is the completion of a unique $\Gm$-equivariant
sheaf of algebras on the product $X^\tw \times \Spec\kk[h]$.
\end{lemma}

\proof{} As in Lemma~\ref{glob}, under the positive weight
assumption, the completion is an equivalence of tensor categories of
$\Gm$-equivariant coherent sheaves on $X^\tw \times \Spec \kk[h]$
and on $\widehat{X}^\tw$.  \endproof

Frobenius-constant quantizations are local objects, both in Zariski
and in \'etale topology. For an affine $X$, the functor of global
sections provides an equivalence between Frobenius-constant
quantizations of $X$ and those of $\calo_X$.

A motivation for the definition is provided by the next

\begin{prop}\label{Dhm} 
If $M$ is a smooth $k$-variety, then $D_h(M)$ is a
Frobenius-constant quantization of $T^*M$.
\end{prop}
 
\proof[Sketch of a proof.] We define a map from generators of
$\calo(T^*M)^p$ to $D_h(M)$. A function $f^p\in \calo_M^p$ lifted
from $M$ is sent to $f^p\in D_h(M)$. A fiber-wise linear function on
$T^*M$ is a vector field on $M$; for such a function $\xi\in
Vect(M)$ we set $\xi^p\mapsto \xi^p-h^{p-1}\xi^{[p]}$, where
$\xi^{[p]}$ is the restricted power of a vector field $\xi$. The
remaining part of the proof is a direct computation in local
coordinates (cf. \cite[\S 1.2]{MR}).
\endproof

\begin{example}\label{Weyl algebra}
The $h$-adic completion $\W$ of the Weyl algebra $\bW$ (see
Subsection~\ref{nota}) is a Frobenius-constant quantization of the
symplectic space $V$. (This may be viewed as a particular case of
Proposition~\ref{Dhm} where $M$ is a Lagrangian subspace in $V$). It
is graded. More precisely, $\bW$ carries a grading such that linear
functions have degree $1$, and $\deg(h)=2$; it induces an action of
the multiplicative group on $\W$, and $\bW$ is identified with the
space of $\Gm$ finite vectors in $\bW$.
\end{example}

\begin{example}\label{Lie alg}
Let $\g$ be a Lie algebra over $\kk$. Then the {\it asymptotic
enveloping algebra} $U_h(\g)$ is the $h$-completion of the graded
algebra $\kk[h]\langle \g\rangle/(xy-yx=h[x,y])$.  If $\g$ is a
restricted Lie algebra, then $U_h(\g)$ is a Frobenius-constant
quantization of the Poisson variety $\g^*$; the Frobenius center is
generated by $x^p-h^{p-1}x^{[p]}$, $x\in \g$.  (Notice that these
central elements are homogeneous with respect to the natural
grading, in which $deg(h)=1$).
\end{example}

The following result is a generaliztion of a fundamental observation
by Mirkovic and Rumynin.

\begin{prop}\label{azu}
Let $\calo_h$ be a Frobenius-constant quantization of a symplectic
variety $X$. Then the restriction $\calo_h[h^{-1}]$ of $\calo_h$ to
the generic fiber of the formal scheme $\widehat X$ is an Azumaya
algebra of rank $p^{\dim X}$.
\end{prop}

\proof{} For a closed point $x \in X$, its Frobenius neighborhood
$Frob(x) \subset X$ is by definition the spectrum of the fiber of
the $\calo_X^p$-algebra $\Fr_*\calo_X$ at the point $x$. It is a
Poisson scheme. It sufficies to check that for any closed point $x
\in X$ and for any quantization $A$ of the Poisson algebra of
functions on $\FrN(x)$, the localization $A(h^{-1})$ is an Azumaya
algebra over $\kk((h))$. Let $Z$ be the center of $A(h^{-1})$. Then
$$
\dim_{\kk((h))} Z =\dim_\kk (Z\cap A) \mod h A.
$$
The right hand side is contained in the Poisson center of
$\calo_{\FrN(x)}$ which is easily seen to be one dimensional. Thus
$Z=\kk((h))$. Assume now that $I\subset A(h^{-1})$ is a non-zero
two-sided ideal. Then $(I\cap A) \mod h$ is a Poisson ideal in
$\calo_{\FrN(x)}$. Since Hamiltonian vector fields generate the
tangent space to $X$ at $x$ (this is immediate by a computation in
local coordinates), it follows that a Poisson ideal is invariant
under all derivations; on the other hand, one checks easily that
$\calo_{\FrN(x)}$ has no proper nonzero ideals invariant under all
derivations. Thus $(I\cap A) \mod h=A\mod h$, so that $I=A$ by
Nakayama lemma.
\endproof

\section{Generalities on symplectic resolutions.}\label{gen}

We return to the set-up of the introduction: $V$ is a
finite-dimensional symplectic vector space, $\Grp \subset Sp(V)$ is
a finite subgroup, $\pi:X \to Y=V/\Grp$ is a smooth projective
resolution, the symplectic form $\omega$ on $X_0 = V_0/\Grp$ extends
to a non-degenerate symplectic form on $X$.

We set $\calo^\Grp=\calo(V)^\Grp=H^0(X,\calo)$. Let the
multiplicative group $\Gm$ act on the vector space $V$ by
dilations. This actions descends to an action on $V/\Grp$, and then
lifts to an action on $X \to V/\Grp$ (see \cite{Ka}; the assumption
$\cchar\kk=0$ adopted in \cite{Ka} is not essential). All the
actions are positive-weight. It is also known \cite{Ka} that the map
$X \to Y$ is necessarily semismall --- that is, $\dim X \times_Y X =
\dim X$ (again, the proof in \cite{Ka} assumes characteristic $0$,
but it works without any changes in arbitrary characteristic).

In addition, we will assume that $X$ satisfies the following:
\begin{itemize}
\item We have $H^p(X,\Omega^q)=0$ when $p+q > \dim X$.
\end{itemize}
This is true in characteristic $0$, see \cite{K3}. It is probably
also true in positive characteristic under some additional
assumptions. We could not find these assumptions, unfortunately;
therefore we simply impose \thetag{$\bullet$} as a standing
assumption on the resolution $X$.

\begin{example}
If $\dim V=2$, then every quotient $V/\Grp$ with $\Grp \subset Sp(V)
= SL(V)$ admits a unique resolution $X \to V/\Grp$ satisfying all
our assumptions. This situation is completely classic; we will give
some more details in Subsection~\ref{dim2}.
\end{example}

Recall that we have dense open subsets $V_0 \subset V_1 \subset
\dots \subset V$, and we denote $X_i = \pi^{-1}(V_i/\Grp)$. We will
need a particular \'etale covering of the open subset $X_1 \subset
X$. Consider the connected components $H_\alpha \subset (V_1/\Grp)$
of the complement $(V_1/\Grp) \setminus (V_0/\Grp)$. Every such
$H_\alpha$ is a closed subscheme in $V_1/\Grp$ of pure codimension
$2$. There exists a linear subspace $V_\alpha \subset V$ of
codimension $2$ such that $H_\alpha=\eta(V_\alpha\cap V_1)$, and the
subgroup $\Grp_\alpha \subset \Grp$ of elements in $\Grp$
stabilizing each vector in $V_\alpha$ is non-trivial.

Denote by $\Grp'_\alpha \subset \Grp$ the stabilizer of the subspace
$V_\alpha$; then $\Grp_\alpha$ is a normal subgroup in
$\Grp'_\alpha$.  Consider $V$ as a $\Grp_\alpha$-module. It splits
as an orthogonal sum $V=W\oplus V_\alpha$, where $\dim W = 2$ (so
that $\Grp_\alpha \subset SL(W)$). The quotient group $N_\alpha =
\Grp_\alpha'/\Grp_\alpha$ acts on $W/\Grp_\alpha$ and on $V_\alpha$;
the second action is a symplectic action on a symplectic vector
space. In keeping with our general notation, let $V_{\alpha,0}
\subset V_\alpha$ be the open subset of vectors with trivial
stabilizer in $N_\alpha$.

The projection $\eta:V \to V/\Grp$ induces a natural map
$$
\eta_\alpha:(W/\Grp_\alpha) \times V_{\alpha,0} \to
\left(\left(W/\Grp_\alpha\right) \times V_{\alpha,0}\right)/N_\alpha
\to V/\Grp.
$$
The map $\eta_\alpha$ is \'etale outside of some closed subset $F_1
\subset (W/\Grp_\alpha) \times V_{\alpha,0}$ which is disjoint from
$\{0\} \times V_{\alpha,0}$.
The preimage $\eta_\alpha^{-1}(H_\alpha)$ is the disjoint union of
$(\{0\} \times V_{\alpha,0})
/N_\alpha$ and some closed subset $F_2 \subset (W/\Grp_\alpha)
\times V_{\alpha,0}$.

Denote by $U_\alpha \subset (W/\Grp_\alpha) \times V_{\alpha,0}$ the
complement to the closed subset $F_1 \cup F_2$. Then $U_\alpha$ is a
Zariski neighborhood of $\{0\} \times V_{\alpha,0} \subset
(W/\Grp_\alpha) \times V_{\alpha,0}$; the map $\eta_\alpha:U_\alpha
\to V/\Grp$ is \'etale, and $\eta_\alpha^{-1}(H_\alpha) = \{0\}
\times V_{\alpha,0} \subset U_\alpha$. Moreover, $\eta_\alpha$
induces an isomorphism $\eta_\alpha:(\{0\} \times
V_{\alpha,0})/N_\alpha \cong H_\alpha$.

The map $\eta_\alpha:(W/\Grp_\alpha) \times V_{\alpha,0} \to V/\Grp$
extends to a map
$$
\rho_\alpha:\left(Y_\alpha \times V_{\alpha,0}\right)/N_\alpha \to
X,
$$
where $\pi_\alpha:Y_\alpha \to W/\Grp_\alpha$ is the canonical
crepant resolution of $W/\Grp_\alpha$ (see e.g. \cite[Section
4]{Ka}; since the resolution $Y_\alpha \to W/G_\alpha$ is
canonical, the action of $N_\alpha$ on $W/\Grp_\alpha$ lifts
uniquely to an action on $Y_\alpha$).

Denote by $X_\alpha \subset Y_\alpha \times V_{\alpha,0}$ the
preimage $(\pi_\alpha \times \id)^{-1}(U_\alpha)$ of the open subset
$U_\alpha \subset (W/\Grp_\alpha) \times V_{\alpha,0}$. The map
$\rho_\alpha$ is \'etale on $X_\alpha$; we have
\begin{equation}\label{un}
\rho_\alpha^{-1}\pi^{-1}(H_\alpha) = \pi_\alpha^{-1}(\{0\}) \times
V_{\alpha,0},
\end{equation}
and the map $\rho_\alpha$ induces an isomorphism
$$
\rho_\alpha:\left(\pi_\alpha^{-1}(\{0\}) \times
V_{\alpha,0}\right)/N_\alpha \to \pi^{-1}(H_\alpha) \subset X.
$$
The sets $X_\alpha$ together with $X_0 \subset X_1$ form an \'etale
covering of the subset $X_1 \subset X$. The following glueing lemma
formalizes the situation.

\begin{lemma}\label{glue}
The category of coherent sheaves $\E$ on $X_1$ is equivalent to the
category of the following data:
\begin{enumerate}
\item a coherent sheaf $\E_0$ on $X_0$,
\item for every component $H_\alpha$ of the complement $V_1/\Grp
\setminus V_0/\Grp$, an extension $\E_\alpha$ of the sheaf
$\rho_\alpha^*\E_0$ to an $N_\alpha$-equivariant sheaf on
$X_\alpha$.
\end{enumerate}
\end{lemma}

\proof{} Since $X_0$ and all the $X_\alpha$ together form an \'etale
covering of $X_1$, it suffices to extend $\E_0$ to a sheaf on
$\rho_\alpha(X_\alpha)$ for every $X_\alpha$. To do this, it
suffices to construct the descent data for the \'etale map
$\rho_\alpha:X_\alpha \to X$. Let $X_{\alpha,0} \subset X_\alpha$ be
the preimage $\rho_\alpha^{-1}(X_0)$. By \eqref{un}, the product
$X_\alpha \times_X X_\alpha$ is covered by two open subsets:
$X_{\alpha,0} \times_X X_{\alpha,0}$ and $X_\alpha \times
N_\alpha$. The descent data on the first subset are already given,
since $\E_\alpha = \rho^{-1}\E_0$ on $X_{\alpha,0} \subset
X_\alpha$. The descent data on the second subset are equivalent to
the $N_\alpha$-equivariant structure on $\E_\alpha$.
\endproof

\begin{prop}\label{lifts}
\begin{enumerate}
\item Assume that Theorem~\ref{evb} holds for some vector bundle
$\E$. Then $\E$ carries a canonical $\Gm$-equivariant structure.
\item Let $R$ be a regular complete local ring with residue field
$\kk$, the maximal ideal $\m \subset R$, and and fraction field
$\K$, and assume given a projective crepant resolution $\pi:X_R \to
V_R/\Grp$, where $X_R$ is a scheme smooth over $R$. Assume that
Theorem~\ref{evb} holds for the special fiber $X_\kk$ of the scheme
$X_R$, and that for every one of the \'etale charts $X_\alpha
\subset Y_\alpha \times V_{\alpha,0}$, the pullback $\rho^*\E$
extends to a vector bundle on the whole $Y_\alpha \times
V_\alpha$. Then Theorem~\ref{evb} holds for the generic fiber
$X_\K$.
\end{enumerate}
\end{prop}

\proof{} To prove \thetag{i}, note that since we have $\calo_V
\subset \calo_V\#\Grp$, any object $\E$ satisfying the conditions of
Theorem~\ref{evb} comes equipped with a structure of an
$\calo_V$-module. In particular, on $X_0 \cong V_0/\Grp \subset X$
we must have
$$
\E \cong \eta_*\wt{\E}
$$
for some sheaf $\wt{\E}$ on $V_0$. Since $\rk \E = |\Grp|$ is equal
to the degree of the \'etale map $\eta:V_0 \to V_0/\Grp$, the sheaf
$\wt{\E}$ must be a line bundle. Since the complement to $V_0
\subset V$ is of codimension $\geq 2$, the Picard group
$\Pic(V_0)=\Pic(V)$ is trivial; therefore $\wt{E}=\calo_V$ and $\E
\cong \eta_*\calo_V$ on $X_0 \subset X$. In particular, it carries a
$\Gm$-equivariant structure on the open part $X_0 \subset X$.

Let $a:\Gm\times X\to X$ be the action map, and let $p:\Gm\times
X\to X$ be the projection. Then an equivariant structure on $\E$ is
given by an isomorphism $a^*\E \cong p^*\E$; we are given such an
isomorphism on $\Gm \times X_0$, and we have to show that it extends
to $\Gm \times X$. We will show that {\em every} isomorphism
$a^*(\E) \cong p^*(\E)$ defined on $\Gm \times X_0$ extends to $\Gm
\times X$. Indeed, let $\X$ be the formal neighborhood of $\{1\}
\times X \subset \Gm \times X$, and let $\X_0$ be the formal
neighborhood of $\{1\} \times X_0 \subset \Gm \times X_0$; then it
suffices to show that every isomorphism $a^*\E \cong p^*\E$ on
$\X_0$ extends to $\X$. Since $\Ext^1(\E,\E)=0$, the sheaf $\E$
extends uniquely to the (trivial) one-parameter deformation $\X$ of
the variety $X$. Therefore there exists at least some isomorphism
$a^*\E \cong p^*\E$. But since $\End(\E) = \calo_V\#\Grp$ is the
same as $\End(\E_{X_0})$, the Formal Function Theorem show that
every section of $\uHom (\E,\E)$ on $\X_0$ extends to $\X$;
therefore any given isomorphism $a^*\E \cong p^*\E$ on $\X_0$ indeed
extends to $\X$.

To prove \thetag{ii}, assume that Theorem~\ref{evb} holds for
$X_\kk$. By \thetag{i} the corresponding vector bundle $\E$ on
$X_\kk$ is $\Gm$-equivariant. Therefore by Lemma~\ref{glob} it
extends to an almost exceptional vector bundle $\E_R$ on $X_R$.

We know that the algebra $\End_R(\E_R)$ is a flat
deformation of the algebra $\End(\E)$. It remains to prove that it
in fact coincides with $Sym(V_R)\#\Grp$.

Applying Lemma~\ref{glob} to $\Gm$-equivariant sheaves on
$V_R$, we see that it suffices to prove that
$$
\End(\E_\X) \cong \wh{Sym(V_R)\#\Grp},
$$
where $\X$ is the formal completion of $X_R$ along $X_\kk \subset
X_R$, and the completion on the right-hand side is taken with
respect to the $\m$-adic topology. Assume first that $V = W \oplus
V'$ with $\dim W = 2$, and that the group $\Grp$ acts trivially on
$V'$, so that, in parituclar, $X = X_1$. Then $\E \cong
\eta_*\calo_V$ on $X_0 = V_0/\Grp \subset X_\kk$, and all the
$\Gm$-homogeneous elements in the group
$$
\Ext^1(\E|_{X_0},\E|_{X_0})
$$
have strictly negative weight with respect to the
$\Gm$-action. Therefore by standard deformation theory, the sheaf
$\E$ on $X_0$ admits a {\em unique} $\Gm$-equivariant extension to a
sheaf on $\X_0$. This means that
$$
\E_{\X_0} \cong \wh{\eta_*\calo_{V}}
$$
on $\X_0 \subset \X$. Now, the restriction to $\X_0 \subset X$
induces a map 
$$
\End(\E_\X) \to \End(\E_{\X_0}) \cong \wh{Sym(V_R)\#\Grp}
$$
of Noetherian topological $R$-algebras complete with respect to the
$\m$-adic topology, and this map is an isomorphism modulo
$\m$. Therefore it is an isomorphism.

In the general case, this argument does not work, because the
extension group $\Ext^1(\E_{X_0},\E_{X_0})$ becomes too
large. However, we can consider the restriction $\E_{X_1}$ of the
bundle $\E$ to $X_1 \subset X_\kk$, and we claim that
$$
\Ext^i(\E_{X_1},\E_{X_1}) = 0
$$
for $i=1,2$. Indeed, consider the vector bundle $\eend(\E)$ on
$X_1$. Theorem~\ref{evb}~\thetag{ii} implies that
$R^k\pi_*\eend(\E)=0$ for $k \geq 1$, so that by
Theorem~\ref{evb}~\thetag{i} we have
$$
\Ext^i(\E,\E) \cong H^i(X_1,\eend(\E)) \cong
H^i(V_1/\Grp,\calo_V\#\Grp).
$$
The algebra $\calo_V\#\Grp$ considered as a sheaf on $V_1/\Grp$ is
the direct image of the trivial sheaf $\calo_V \otimes \kk[\Grp]$
with respect to the quotient map $\eta:V_1 \to V_1/\Grp$. Therefore
$$
H^i(V_1/\Grp,\calo_V\#\Grp) = H^i(V_1,\calo_V) \otimes \kk[\Grp].
$$
Since $\codim V \setminus V_1 \geq 4$, the right-hand side is indeed
trivial for $i = 1,2$.

Now, by the same argument as in the case $V=W \oplus V'$, it
suffices to prove that $\End(\E_{\X_1})$ coincides with the
smash-product algebra $\wh{Sym(V_R)\#\Grp}$. Since
$\Ext^1(\E_{X_1},\E_{X_1}) = 0$, the vector bundle $\E$ on $X_1$
admits a unique extension to a vector bundle on $\X_1$. Thus it
suffices to construct at least {\em some} vector bundle $\E'$ on
$\X_1$ which extends $\E_{X_1}$ and has the correct endomorphism
algebra. To do this, we apply Lemma~\ref{glue}, and deduce that it
suffices to construct a $N_\alpha$-equivariant extension $\E'_\alpha$
of the sheaf $\E_\alpha$ on each of the resolution $Y_\alpha \times
V_\alpha \to V/\Grp_\alpha$ corresponding to the \'etale charts
$X_\alpha$. But since $\Grp_\alpha \subset \Grp$ acts trivially on
$V_\alpha$, we have already proved that the extension $\E'_\alpha$
exists, and it is unique. In particular, it is
$N_\alpha$-equivariant.
\endproof

\section{Quantization of a symplectic resolution}\label{quant}

Assume now that the base field $\kk$ is a perfect field of
characteristic $p > 0$; we will also assume that $p > \dim X$, $p>
|\Grp|$. The goal of this section is the following

\begin{theorem}\label{thmqua} 
In the above assumptions there exists a graded
Frobenius-con\-s\-tant quantization $\calo_h$ of $X$ such that the
$(\calo^\Grp)^p$-algebra $H^0(X,\calo_h)$ of global section of
$\calo_h$ is isomorphic to the standard quantization $\W^\Grp$;
the isomorphism is compatible with the $\Gm$-actions.
\end{theorem}

\begin{lemma}\label{H^0} 
Assume that a graded Frobenius-constant quantization $\calo_h$ of
$X$ is such that its restriction to the open stratum $X_0\cong
V_0/\Grp$ is $\Gm$-equivariantly isomorphic to
$\W|_{V_0^\tw}^\Grp$. Then $\calo_h$ satisfies the conditions of
Theorem~\ref{thmqua}.
\end{lemma}

\proof{} Since the complement to $V_0$ has codimension at least 2,
we have $\W=H^0(\W|_{V_0^\tw})$. Thus $H^0(\calo_h)$ is a subalgebra
in $\W^\Grp=H^0(X_0,\calo_h|_{X_0^\tw})$. Using vanishing of
$H^1(X,\calo)$ we see that
$$
H^0(\calo_h)/h H^0(\calo_h)\iso H^0(\calo)
=\calo_V^\Grp=\W^\Grp/h\W^\Grp,
$$
which shows surjectivity of the map $H^0(\calo_h)\to \W^\Grp$.
\endproof

The plan of the proof of the Theorem is as follows. By
Lemma~\ref{H^0}, it suffices to take the given quantization of $X_0$
and extend it to a quantization of $X$. In Subsection~\ref{toX} we
will extend this quantization to $X_1$, and then to the whole of
$X$. Since $\codim (V_1-V_0)= 2$, extension to $X_1$ reduces to the
case $\dim (V)=2$ which is treated separately in
Subsection~\ref{dim2}.

\subsection{Dimension $2$ case.}\label{dim2}

Throughout this subsection we assume that the dimension $\dim V =
2$. Our aim here is

\begin{prop}\label{dim2 Prop}
Theorem~\ref{thmqua} holds if $\dim (V)=2$. Moreover, the resolution
$X$ and quantization $\calo_h$ are equivariant under the normalizer
of $\Grp$ in $Sp(V)$.
\end{prop}

\subsubsection{The resolution as a Hamiltonian reduction.}

Classification of the data $\langle V,\omega,\Grp\rangle$ with
$\dim(V)=2$ is well-known (and goes back at least to Klein). They
are in bijection with simply-laced Dynkin diagrams; the crepant
resolution $X$ of $V/\Grp$ is unique in each case. Recall a recent
description of $X$ as the Hilbert quotient of $V$ by $\Grp$ due to
Ito and Nakamura \cite{IN} (see also the exposition in
\cite{N1}). Consider the subvariety of $\Grp$-fixed points in the
Hilbert scheme $\Hilb^n(V)$ of subschemes in $V$ of length $n$,
$n=|\Grp|$, and let $\Hilb^\Grp(V)$ be its connected component which
contains elements of the form $\bigoplus_{\gamma \in \Grp}
\calo_{\gamma(v)}$, $v\in V\setminus \{0\}$. Then we have $X\cong
\Hilb^\Grp(V)$.

The resolution $X$ can also be obtained as a Hamiltonian reduction
of an (algebraic) symplectic vector space by an action of a
reductive group; this construction (in analytic set-up) was
discovered by Kronheimer \cite{Kr}, and studied in works of Lusztig,
Nakajima and others. Let us briefly describe a version of this
construction (apparently due to Nakajima) based on the realization
as the Hilbert quotient.

Let $R=\kk[\Grp]$ be the regular representation of the group
$\Grp$. The space $X=Hilb^\Grp(V)$ is the moduli space of
representations of $Sym(V^*)\#\Grp$ which are isomorphic to $R$ as a
$\Grp$-module, and are generated by a $\Grp$-invariant vector.  Let
$\R$ denote the tautological vector bundle on $X$, i.e. the
pushforward to $X$ of the structure sheaf of the universal
$\Grp$-equivariant subscheme in $X\times V$.

Let $M$ be the space parametrizing representations $N$ of
$T(V^*)\#\Grp$ equipped with an isomorphism of $\Grp$-modules
$N\cong R$; here $T$ denotes the tensor algebra. By definition, $M$
is the vector space
$$
M = \Hom_\Grp(V^*,\End(R)) = \Hom_\Grp(V^* \otimes R, R) = (V
\otimes \End(R))^\Grp
$$
of $\Grp$-equivariant maps from $V^*$ to $\End(R)$.

(One can also describe $M$ as the space of representations of the
so-called ``doubled quiver'' corresponding to the affine Dynkin
diagram of type $A$, $D$, $E$ respectively; see e.g. \cite{Lu}).

The vector space $\End(R)$ is equipped with the symmetric trace
pairing $\tr(ab)$ and a Lie bracket $[-,-]: \End(R) \otimes
\End(R)\to \End(R)$; both are $\Grp$-invariant.  The pairing
tensored with the symplectic form $\omega$ on $V$ gives a symplectic
form $\Omega$ on $V \otimes \End(R) $; restricting it to the space
$M$ of $\Grp$ invariants we get a symplectic form $\Omega$ on
$M$. The bracket tensored with the form $\omega$ gives a quadratic
map $V \otimes \End(R)\to \End(R)$; let $\tilde\mu:M \to \End(R)$
denote its restriction to $M$.

The group $G=\Aut_\Grp(R)/\Gm$ acts on $M$ preserving $\Omega$.
Let $\g$ be the Lie algebra of $G$. It is easy to check that the
image of $\tilde \mu$ is contained in the subspace $\g^* \subset
(\End(R))^*= \End(R)$; thus we have a map $\mu:M\to \g^*$. One
readily checks that $\mu$ is the moment map for the action of $G$ on
$M$.

The zero fiber $\mu^{-1}(0)$ parametrizes those representations of
$T(V^*)\#\Grp$ which factor through $Sym(V^*)\#\Grp$; thus a point
of $\mu^{-1}(0)$ defines a $\Grp$-equivariant coherent sheaf on
$V$. It turns out that one of the GIT quotients of $\mu^{-1}(0)$ by
$G$ coincides with $X$. Namely, consider the splitting $\iota:G\to
\Aut(R)$ of the projection $\Aut(R) \to G$ which identifies $G$ with
the subgroup of automorphisms which are trivial on the
$1$-dimensional space of $\Grp$-invariant vectors. Define a
character $\chi:G\to \Gm$ by $\chi(g)=\det(\iota (g))$. Then a point
$x\in \mu^{-1}(0)$ is $\chi$-stable iff the corresponding
$\Grp$-equivariant coherent sheaf on $V$ is generated by a
$\Grp$-invariant section (see e.g. \cite{Na98}, cf. also
\cite{King}). Thus the symplectic variety $X=\Hilb^\Grp(V) =
\mu^{-1}(0) \sllash{}_\chi G$ is the Hamiltonian reduction of $M$ by
$G$.

In these terms one can also describe the tautological vector bundle
$\R$. To this end, consider the trivial vector bundle on $M$ with
fiber $R$. The splitting $\iota$ can be used to equip it with a
$G$-equivariant structure: we define an action of $G$ on $R$ by
pulling-back the tautological action of $G$ on $R$ via $\iota$; we
then let $G$ act on $R\otimes \calo$ diagonally.

The restriction of the resulting $G$-equivariant bundle to
$\mu^{-1}(0)$ descends to a vector bundle on the quotient. The
resulting vector bundle on $X$ is identified with $\R$.

As usual in the Hamiltonian reduction picture, we can also do the
steps in reverse order: first take the GIT quotient
$M\sllash{}_{\chi}G$, then $X$ is realized as a closed subscheme in
$ M\sllash{}_{\chi}G$.

To sum up, we have

\begin{prop}\label{pseudoKing}
\begin{enumerate}
\item 
Let $M^{ss}\subset M$ be the open subscheme of representations
generated by a $\Grp$-invariant vector. Then the action of $G$ on
$M^{ss}$ is free, and the geometric quotient $M^{ss}/G$ exists.
%, the projection $M^{ss} \to M/G$
%  identifies $M^{ss}$ with the total space of a locally trivial
%  principal $G$-bundle over $M/G$, and 
%Vrode kak v Mamforde napisano, chto eti 2 usloviya
% i znachat, chto proekciya v faktor
%-- rassloenie.
The subvariety
$$
\left(\mu^{-1}(0) \cap M^{ss}\right)/G \subset M^{ss}/G
$$
is isomorphic to $X$, and the reduction of the form $\Omega$ equals
$\omega$.
  
\item Equip the trivial vector bundle $R\otimes \calo$ with the
  diagonal $G$-action. Then the descent of $R\otimes \calo_{M^{ss}}$
  to $M^{ss}/G$ restricts to $\R$ on $X \subset M^{ss}/G$.\endproof
\end{enumerate}
\end{prop}

\subsubsection{Quantum version.}

To obtain quantizations by Hamiltonian reduction, we introduce the
following

\begin{defn}
  Let $X$ be a symplectic algebraic variety over a field $\kk$ of
  positive characteristic. Let $G$ be an algebraic group acting on
  $X$. A Frobenius-constant quantization $\calo_h$ of $X$ is {\em
  Frobenius $G$-constant} if $G$ acts on $\calo_h$ so that the
  action of $G$ on the central subalgebra $\calo^p_X[[h]]\subset
  \calo_h$ fixes $h$ and coincides with the natural $G$-action on
  $\calo^p_X \subset \calo^p_X[[h]]$.
\end{defn}

We will abuse terminology by saying ``a map from a vector space $W$
to a sheaf $\F$'' instead of ``a map from $W$ to the space of global
sections of $\F$''.

\begin{defn}
  Let $X$, $G$ be as above, let $\g$ be the Lie algebra of the group
  $G$, and let $\calo_h(X)$ be a Frobenius $G$-constant quantization
  of $X$.  A {\it quantum moment map} is a $\kk[[h]]$-algebra
  homomorphism $\mu:U_h(\g)\to \calo_h(X)$ such that
\begin{enumerate}
\item The restriction of $\mu$ to the central subalgebra
  $Sym(\g^\tw)[[h]]$ maps the subalgebra $Sym(\g^\tw) \subset
  Sym(\g^\tw)[[h]]$ into $\calo(X^\tw) \subset \calo(X^\tw)[[h]]$.
\item For every $\xi\in \g$ and every local section $s$ of $\calo_h$
  we have $\mu(\xi)s-s\mu(\xi)=h \xi(s)$, where $\xi(s)$ is the
  action of the Lie algebra $\g$ on $\calo_h$ induced by the action
  of $G$.
\end{enumerate}
\end{defn}
One checks immediately that by virtue of \thetag{ii}, the induced
map
$$
\mu_0:Sym(\g^\tw) \to \calo_X^p
$$
whose existense is guaranteed by \thetag{i} is a moment map for
the $G$-action on $X$.

\begin{example}\label{linexample}
  Let $X=V$ be a symplectic vector space, and let $G=Sp(V)$. Then we
  have a quantum moment map sending $x\in \ssp(V)\cong
  Sym^2(V)\subset V\otimes V$ into the corresponding element of
  $\W$.
  
  If $G\subset Sp(V)$ is an algebraic subgroup, then we get a
  quantum moment map $\mu:\g\to \W$ by restricting the above map to
  the Lie subalgebra $\g \subset \ssp(V)$.

Notice that the induced homomorphism $U_h(\g)\to \W$ (which we also
denote by $\mu$) is related to the $\Gm$-action by the formula $\mu
(tx)=t^2\mu(x)$.
\end{example}

\begin{prop}\label{Ham_red_char_p}
Let $M$ be a smooth symplectic manifold, let $G$ be a group acting
on $M$, let $\calo_h$ be a Frobenius $G$-constant quantization of
$X$, and assume given a quantum moment map $\mu:\g\to
\calo_h$. Moreover, assume given a $G$-invariant open subset $U
\subset X$ such that the action of $G$ on $U$ is free, and the
geometric quotient $U/G$ exists; thus the projection $\rho:U \to
U/G$ identifies $U$ with the total space of a principal $G$-bundle
on $U/G$.
% Mumford, prop-n 0.9 
Then the sheaf
$$
(\rho_*\calo_h)^G \subset \rho_*\calo_h
$$
is a Frobenius-constant quantization of the Poisson manifold $U/G$,
and its quotient
$$
\calo_h(Z) = (\rho_*\calo_h)^G/\langle
\rho_*\calo_h\mu(\g)\rho_*\calo_h\rangle^G
$$
is a Frobenius-constant quantization of the Hamiltonian quotient
$$
Z = \left(\mu_0^{-1}(0) \cap U\right)/G \subset U/G.
$$
Moreover, if $\F$ is a locally projective $G$-equivariant sheaf of
left $\calo_h$-modules on $X$, then the sheaf
$(\rho_*\F)^G/I(\rho_*\F)^G$ is a locally projective sheaf of
$\calo_h(Z)$-modules on $Z$.\endproof
\end{prop}

\proof[Sketch of a proof.]  If $S\subset U/G$ is an affine open
subset, then, since the action of $G$ on $\rho^{-1}(S)$ is free, the
algebra of functions on $\rho^{-1}(S)$ is an injective object in the
category of algebraic $G$-modules. This proves the first
statement. On the other hand, since the action is free,
$\calo(\rho^{-1}(S))$ is flat over $Sym(\g)$, and $\mu(\g)
\calo(\rho^{-1}(S))$ is an injective $G$-module.  This implies an
isomorphism
\begin{equation}\label{vse_ravno}
\calo_Z=\rho_*(\calo)^G/(\mu(\g)\rho_*(\calo)^G)\iso
\left( \rho_*(\calo)/ \mu(\g)\rho_*(\calo) \right) ^G.
\end{equation}
Also, $\calo_h(\rho^{-1}(S))$ is flat over $U_h(\g)$, hence $\left(
\rho_*(\calo_h)/ \mu(\g)\rho_*(\calo_h) \right) ^G$ is an $h$-flat
deformation of $\calo (Z)$.  The isomorphism \eqref{vse_ravno}
implies
$$
\calo_h(Z)=\rho_*(\calo_h)^G/(\mu(\g)\rho_*(\calo_h)^G) \iso \left(
\rho_*(\calo_h)/ \mu(\g)\rho_*(\calo_h) \right) ^G.
$$
As noted above, the right hand side here is an $h$-flat deformation
of $\calo(Z)$, and the left hand side carries associative
multiplication. It is clear that the commutator map in the resulting
sheaf of algebras is compatible with the Poisson bracket; thus
$\calo_h(Z)$ is indeed a quantization of $\calo_Z$. The proof of the
statement about the Hamiltonian reduction of a $\calo_h$-module is
parallel to the proof of the statement about the reduction of the
sheaf of algebras $\calo_h$.
\endproof

\subsubsection{Quantizing the resolution.}

We now return to the situation of Proposition~\ref{dim2 Prop}.
Combining Proposition~\ref{pseudoKing} with
Proposition~\ref{Ham_red_char_p}, we obtain a
Fro\-be\-ni\-us-con\-s\-tant quantization $\calo_h$ of the
resolution $X \to V/G$. It is clear from the construction that the
quantization $\calo_h$ is graded.  By Lemma~\ref{H^0}, to prove
Proposition~\ref{dim2 Prop}, it remains to show that $\calo_h$
coincides with the standard quantization on $V_0 \subset V$.

Consider the sheaf $\R_h$ of $\calo_h$-modules on $X^\tw$ obtained
as the Hamitonian reduction of the free module $R\otimes
\calo_h(M)$. Then by Proposition~\ref{Ham_red_char_p} together with
Proposition~\ref{pseudoKing}~\thetag{ii}, we have $\R_h/h\R_h\cong
\R$.

Moreover, consider the restriction of the sheaf $\R_h$ to the open
subset $X_0 = V_0/\Grp \subset X$. Then we claim that there exists a
natural isomorphism
$$
\R_h \cong \calo_h \otimes_{\calo_{X_0}^p} \calo_{V_0}^p
$$
of sheaves on $X_0$. Indeed, the definition of $\R_h$ yields an
action of $\calo(V^\tw)[[h]]=H^0(X^\tw,\R^\tw)[[h]]$ on $\R_h$
commuting with the action of $\calo_h$. This action, in turn, yields
a map of sheaves
$$
\calo^p_{V_0}\otimes_{\calo^p_{V/\Grp}} \calo_h\to \R_h.
$$
It is readily seen to induce an isomorphism modulo $h$; since both
sheaves are $h$-flat it is an isomorphism.

To prove Proposition~\ref{dim2 Prop}, it remains to apply the
following

\begin{lemma}
Assume that $\calo_h$ is a $\Gm$-equivariant Frobenius-constant
quantization of $X$ such that the sheaf
$$
\calo_h|_{X_0}\otimes _{\calo_{X_0}^\tw}\calo_{V_0^\tw}
$$
of $\calo_h$-modules on $X_0=V_0/\Grp$ extends to a locally
projective $\calo_h$-module $\R_h$ on $X$ satisfying
$\R_h/h\R_h\cong \R$.

Then the restriction of the quantization $\calo_h$ to the open
subset $X_0^\tw=V_0^\tw/\Grp$ is isomorphic to
$\W|_{V_0^\tw}^\Grp$; the isomorphism is compatible with the
$\Gm$-actions.
\end{lemma}

\proof{} It sufficies to construct a $\Grp$-equivariant isomorphism
of graded $(\calo^\Grp)^p$-algebras
\begin{equation}\label{H0 is W}
H^0\left(\calo_h\otimes_{\calo_{X_0}^\tw}\calo_{V_0^\tw}\right)\cong
\W.
\end{equation}
Consider the subspace $V^* \subset \calo(V)$ of linear functions on
$V$. It follows e.g. from \cite{KV} that $H^i(X,\R)=0$ for $i>0$.
Hence the map
$$
H^0(\R_h)=H^0(\calo_h\otimes_{\calo_{X_0}^\tw}\calo_{V_0^\tw})
\to H^0(\R) =\calo(V)
$$
is surjective. In particular, we can factor the embedding $V
\subset \calo(V)$ through a map $\iota:V^* \to H^0(\R_h)$.
Moreover, we can choose a $\Gm$-equivariant map, so that
$\Im(\iota)$ is contained in the space of elements of degree $1$
with respect to the $\Gm$-action. Then for any $x,y\in V^*$, the
element $xy-yx$ equals $h \omega(x,y)$ modulo $h^2$, and has
degree $2$ with respect to the $\Gm$-action. This clearly implies
$xy-yx =h \omega(x,y)$. Thus we get a $\Grp$-equivariant algebra
homomorphism from the completed Weyl algebra $\W$ to
$$
H^0(\calo_h\otimes_{\calo_{X_0}^\tw}\calo_{V_0^\tw}).
$$
Since it induces an isomorphism on the associated graded spaces with
respect to the $h$-adic filtrations, it is itself an isomorphism.
\endproof

\subsection{The general case.}\label{toX}

Return to the situation of Theorem~\ref{thmqua}, the general case
(no assumption on $\dim V$).  Recall notations $V_1$, $X_1$ from
Subsection~\ref{nota}.

\begin{prop}\label{Prop_on_X1}
The standard quantization $\W^\Grp$ of $X_0 = V_0/\Grp$ extends
to a graded Frobenius-constant quantization $\calo_h$ of $X_1$.
\end{prop}

\proof{} A Frobenius-constant quantization of $X$ is a coherent
sheaf on $X\tw$, so that Lemma~\ref{glue} applies. Thus it suffices
to extends the given quantization of $X_0 = V_0/\Grp \subset X$ to a
quantization of each of the \'etale charts
$X_\alpha$. Propositon~\ref{dim2 Prop} applied to $Y_\alpha$ gives a
graded Frobenius-constant $N_\alpha$-equivariant quantization of
$Y_\alpha$. Taking its ($h$-adically completed) tensor product with
the standard quantization of the vector space $V_\alpha$, we obtain
by restriction a graded Frobenius-constant $N_\alpha$-equivariant
quantization of $X_\alpha \subset Y_\alpha \times V_\alpha$.
\endproof

To finish the proof of Theorem~\ref{thmqua}, it remains to show that
the Frobenius-constant quantization $\calo_h$ of $X_1 \subset X$
extends to a Frobenius-constant quantization of the whole $X$. It
suffices to prove the following.

\begin{prop}\label{allX}
Every graded Frobenius-constant quantization $\calo_h$ of $X_1
\subset X$ extends uniquely to a graded Frobenius-constant
quantization of $X$.
\end{prop}

\proof{} First of all observe that it is enough to extend $\calo_h$
to a deformation of $\calo(X)$ to a sheaf of associative algebras:
then the action of $\Gm$, and the embedding of $\calo^p$ into the
center of $\calo_h$ extend uniquely from $X_1$, because $\codim
(X\setminus X_1) \geq 2$.

We will use Lemma~\ref{def}. The Hochschild cohomology of smooth
varieties can be computed by the Hochschild-Kostant-Rosenberg
spectral sequence, whose $E_2$-term is given by
$$
H^p(X,\Lambda^q\T(X)),
$$
the cohomology of degree $p$ with coefficients in the bundle
$\Lambda^q\T(X)$ of polyvector fields of degree $q$ on $X$. Then
Lemma~\ref{def} reduces the Proposition to the following
cohomological statement.

\begin{lemma}
Consider the restriction map
$$
\sigma:H^p(X,\Lambda^q\T(X)) \to H^p(X_1,\Lambda^q\T(X_1))
$$
on the components of negative weight. Then it is an isomorphism if
$p+q=1,2$, and it is injective if $p+q=3$.
\end{lemma}

\proof{} Since the sheaves of polyvector fields are locally free, while
the complement $X \setminus X_1 \subset X$ is of codimension at
least $2$, the map $\sigma$ is a bijection for $p=0$ and an
injection for $p=1$ (for all weights).

By our assumption \thetag{$\bullet$}, we have
$$
H^p(X,\Omega^q_X) = 0
$$
when $p+q > \dim X$. Since $X$ is symplectic, we have $\Omega^p_X
\cong \Lambda^{\dim X - p}\T(X)$, so that this implies
$$
H^p(X,\Lambda^q\T(X)) = 0
$$
when $p > q$. For $q = 0$ this gives $H^p(X,\calo_X) = 0$, $p \geq
1$, and $H^p(X_1,\calo_X) = H^p(V_1/\Grp,\calo_{V/\Grp}) =
H^p(V_1,\calo_V)^\Grp$. Since $\codim V \setminus V_1 \geq 4$,
this implies that $H^p(X_1,\calo_X)=0$ for $p=1,2$. For $q=1$ we
obtain $H^2(X,\T(X)) = 0$, so that the map $\sigma$ is
tautologically injective for $p=2$, $q=1$.

It remains to consider the sum of components of negative weight in
the group $H^1(X,\T(X))$. We can identify $\T(X) \cong \Omega^1_X$
by means of the symplectic form; since the symplectic form has
weight $2$, we are reduced to studying the sum of components of
non-positive weights in $H^1(X,\Omega^1_X)$. Moreover, since the
$\Gm$-action on $X$ is positive-weight, we only have to analyze the
$\Gm$-invariant part in $H^1(X,\Omega^1_X)$. We already know that
the map $\sigma$ is injective on $H^1(X,\Omega^1_X)$. Since $\codim
X \setminus X_1 \geq 2$, every line bundle on $X_1$ extends to the
whole of $X$. Therefore to finish the proof, it suffices to prove
the following.

\begin{lemma}
  The group $H^1(X_1,\Omega^1_X)^\Gm$ is generated by Chern classes of
  line bundles on $X_1$.
\end{lemma}

\proof{} The sheaf $\pi_*\Omega^1_X$ on $V_1/\Grp$ coincides with
the $\Grp$-invariant part in the sheaf $\eta_*\Omega^1_V$, where
$\eta:V \to V/\Grp$ is the quotient map. Since $\codim V \setminus
V_1 \geq 4$, we have $H^i(V_1/\Grp,\pi_*\Omega^1(X)) = 0$ for $i =
1,2$. Therefore
$$
H^1(X_1,\Omega^1_X) = H^0(V_1/\Grp,R^1\pi_*\Omega^1_X).
$$
The sheaf $R^1\pi_*\Omega^1_X$ is concentrated on the complement
$V_1/\Grp\setminus V_0/\Grp$, which is a disjoint union of the
components $H_\alpha$. Therefore on $V_1/\Grp$ we have
$$
R^1\pi_*\Omega^1_X = \bigoplus_\alpha
\left(\eta_{\alpha,*}\eta_\alpha^* R_1\pi_*\Omega^1_X
\right)^{N_\alpha},
$$
where $\eta_\alpha:U_\alpha \to V_1/\Grp$ are the \'etale charts
constructed in Section~\ref{gen}. By base change we have
$$
\eta_\alpha^*R^1\pi_*\Omega^1_X \cong
R^1\pi_*\Omega^1_{X_\alpha}.
$$
Recall that $X_\alpha$ is an open dense subset in
$\overline{X}_\alpha = Y_\alpha \times V_\alpha$. The sheaf
$\Omega^1_{\overline{X}_\alpha}$ on $\overline{X}_\alpha$ decomposes
as a direct sum
$$
\left( \Omega^1_{Y_\alpha} \boxtimes \calo_{V_\alpha} \right) \oplus
\left( \calo_{Y_\alpha} \boxtimes \Omega^1_{V_\alpha} \right),
$$
and since $H^1(Y_\alpha,\calo_{Y_\alpha})=0$, we have
$$
R^1(\pi_\alpha\times id)_* (\calo_{Y_\alpha} \boxtimes
\Omega^1_{V_\alpha})=0
$$
by the projection formula. Thus 
$$
R^1\pi_*\Omega^1_{X_\alpha} \cong
R^1\pi_{\alpha,*}\Omega^1_{Y_\alpha} \boxtimes \calo_{V_\alpha}.
$$
We conclude that
$$
H^0(V_1/\Grp, R^1\pi_*\Omega^1_{X_1}) \cong
\bigoplus_\alpha \left(H^1(Y_\alpha,\Omega^1_{Y_\alpha}) \otimes
H^0(U_\alpha,\calo_{V_{\alpha,0}})\right)^{N_\alpha},
$$
and since the complement $V_\alpha \setminus V_{\alpha,0}$ is of
codimension at least $2$, we have
$$
H^0(U_\alpha,\calo_{V_{\alpha,0}}) \cong \calo_{V_\alpha}
$$
for each $\alpha$. It remains to prove that the group of weight-$0$
elements in
\begin{equation}\label{grpp}
\left(H^1(Y_\alpha,\Omega^1_{Y_\alpha}) \otimes
\calo_{V_\alpha}\right)^{N_\alpha}
\end{equation}
is generated by Chern classes of $N_\alpha$-equivariant line bundles
on $\overline{X}_\alpha$. But it is easy to see that
$H^1(Y_\alpha,\Omega^1_{Y_\alpha})$ is generated by the Chern
classes of the exceptional curve $E = \pi_\alpha^{-1}(0) \subset
Y_\alpha$, and $\Gm$ acts trivially on this group. On the other
hand, the only functions of non-positive weight on $V_\alpha$ are
constants. Thus every $N_\alpha$-invariant element of weight $0$ in
the group \eqref{grpp} is a linear combination of Chern classes of
divisors of the form $D \times V$, where $D$ is a
$N_\alpha$-invariant divisor supported on $E$. Every such divisor $D
\times V$ corresponds to a divisor on $X$.
\endproof

\section{Equivalences}\label{last}

From now on, subindices will be used to denote the ring of scalars;
thus for an $R$-algebra $A=A_R$, the ring $A\otimes _R R'$ will be
denoted by $A_{R'}$.

\subsection{Assumptions.}\label{ass}

Fix an algebraically closed field $\K$ of characteristic zero, and
the set of data $\langle V,\omega, \Grp, X\rangle$ defined over $\K$
-- $V$ is a finite-dimesnional vector space over $\K$, $\omega$ is a
symplectic form on $V$, $\Grp \subset Sp(V)$ is a finite subgroup,
and $X \to V/\Grp$ is a projective smooth crepant resolution of the
quotient $V/\Grp$. Let $R\subset K$ be a $\Zet$-algebra of finite
type such that $\langle V,\omega, \Grp, X\rangle$ are defined over
$R$; we assume $R\cong {\mathbb A}^n_R$, $\omega_R$ is symplectic;
$\Grp\subset Sp(n,R)$, $X$ is smooth over $R$, and $\pi: X\to
V/\Grp$ is proper (below we will just say ``defined over $R$''
leaving the natural conditions implicit). Then the Weyl algebra
$\bW$ is also defined over $R$.  It is clear that such an $R$
exists.

By \cite{Ka} the map $X \to V/\Grp$ is semismall. By \cite{K3}, this
implies that $H^p(X,\Omega^q_X) = 0$ when $p + q > \dim
X$. Therefore localizing $R$ if necessary, we may assume that the
reduction $X_\kk$ to any closed point $\Spec \kk \in \Spec R$ of
positive characteristic satisfies the assumption \thetag{$\bullet$}
in the beginning of Section~\ref{quant} (to check that the bad loci
are Zariski-closed algebraic subvarieties, one uses the
positive-weight $\Gm$-action). Localizing $R$ even further, we can
insure that in fact all the assumptions of Section~\ref{quant} are
satisfied for every $X_\kk$.

Denote $\bWbar =\bW/(h-1)\bW$. 

\begin{lemma}\label{Morita} 
There exists a dense affine open subset $\Spec(R')\subset \Spec(R)$
such that the bimodule $\bWbar_{R'}$ gives a Morita equivalence
$$ 
\bWbar^\Grp_{R'} \sim (\bWbar\#\Grp)_{R'}.
$$
\end{lemma}

\proof{} We can assume that $|\Grp|^{-1}\in R$. Set
$$
e=\frac{1}{|\Grp|}\sum\limits_{\grp\in \Grp}\grp \in R[\Grp]\subset
(\bWbar\#\Grp).
$$
Then $e$ is an idempotent, $\bWbar^\Grp_{R'} =e\left(
(\bWbar\#\Grp)_ {R'}\right) e$, and it is well-known that it
suffcies to check that $e$ generates the unit two-sided ideal in
$(\bWbar\#\Grp)_{R'}$. The algebra $\bWbar_{\K}$ is simple; it is
easy to deduce that the algebra $(\bWbar_{\K}\#\Grp)$ is also
simple; so that $e$ generates the unit ideal in
$(\bWbar\#\Grp)_{\K}$. Hence it also generates the unit ideal in
$(\bWbar\# \Grp)_{R'}$ for some $R'\subset K$ of finite type over
$R$.  \endproof

\subsection{Equivalence for the Azumaya algebra in positive characteristic.}

Take $R'$ satisfying Lemma~\ref{Morita}, and fix a closed point
$\Spec\kk \in \Spec R'$. Let $X = X_\kk$, $V = V_\kk$. Consider the
graded Frobenius-constant quantization $\calo_h$ of $X$ constructed
in Theorem~\ref{thmqua}. Extend it to a sheaf $\bO_h$ on $X \times
\Spec \kk[h,h^{-1}]$ by Lemma~\ref{gra}. Set
$\bO=\bO_h/(h-1)\bO_h$; thus, $\bO$ is a locally free sheaf on
$X^\tw$.

\begin{lemma}
We have $H^i(\bO)=0$ for $i>0$, and $H^0(\bO)=\bWbar^\Grp$ as
a $\calo(V^\tw)^\Grp$-algebra.
\end{lemma}

\proof{} This statement immediately follows from the constructionl
however, it deserves to be stated as a Lemma, since it is the
crucial point in the proof of Theorem~\ref{main}.
\endproof

By Proposition~\ref{azu}, the algebra $\bO$ is an Azumaya algebra
of rank $p^{\dim X}$.

\begin{theorem}\label{equivAzu}
We have equivalences of derived categories
$$
D^b\left(\Coh(X,\bO)\right) \cong
D^b\left(\bWbar_{\kk}^\Grp\fmod\right)\cong
D^b\left(\bWbar{\#\Grp}_{\kk}\fmod\right),
$$
where the first equivalence is given by the derived functor of
global sections, and the second one is provided by
Lemma~\ref{Morita}.
\end{theorem}

\proof{} The second equivalence is given by Lemma~\ref{Morita}.  It
remains to show that the derived functor of global sections provides
an equivalence
$$
D^b\left(\Coh(X,\bO)\right) \cong
D^b\left(\bWbar_{\kk}^\Grp\fmod\right).
$$
We will check that the free rank $1$ module over $\bO$ is an
almost exceptional object in $\Coh(X^\tw,\bO)$; then the desired
statement follows from Proposition~\ref{CYProp}.

The algebra $\bWbar$ has finite homological dimension, since its
associated graded algebra with respect to the natural filtration is
the symmetric algebra $Sym(V^*)$, which has finite homological
dimension. Since $\cchar\kk$ does not divide $|\Grp|$, the algebra
$\bWbar\#\Grp$ also has finite homological dimension; by the
equivalence of Lemma~\ref{Morita}, the same is true for the algebra
$\bWbar^\Grp$. Since we have $H^i(X,\calo)$ for $i \geq 1$, the
object $\bO\in \Coh(X^\tw_{\kk },\bO)$ is indeed almost exceptional.
\endproof

\subsection{Untwisting.} Keep the notation of the last
Subsection. We will now eliminate the Azumaya algebras.

\begin{prop}\label{untwisting}
There exists an Azumaya algebra $\A$ on $V^\tw/\Grp$ whose
restriction to $V_{0}^\tw/\Grp$ is equivalent to
$W^\Grp|_{V^\tw_0/\Grp}$.
\end{prop}

\noindent
The Proposition follows directly from the next

\begin{lemma}
Let $S$ be an affine scheme acted upon by a finite group $\Grp$
with $|\Grp|$ invertible on $S$. Let $p:S\to S/\Grp$ be the projection
to the categorical quotient. Let $l\in \Zet_{>0}$ be prime to $|\Grp|$.
Then $p^*$ induces an isomorphism $Br(S/\Grp)[l]\iso Br(S)[l]^\Grp$
where $[l]$ denotes the $l$-torsion.
\end{lemma}

\proof{} By a Theorem of O. Gabber \cite{Ga}, for any affine scheme
$X$ we have $Br(X)[l]=H^2(X_{et},\Gm)[l]$. For any finite morhism
$p$, the higher direct images in the \'etale topology $R^ip_*(\Gm)$,
$i \geq 1$ are trivial; indeed, the stalk of $R^ip_*(\Gm)$ at a
geometric point $x$ is the $i$-th cohomology group of the \'etale
sheaf $\calo^*$ on the spectrum of some commutative ring $R$ which
is finite over a strictly Henselian ring (see \cite[Theorem
III.1.15]{Milne}); this cohomology group vanishes since $R$ is
itself strictly Henselian by \cite[Corollary I.4.3]{Milne}.
 
We have the norm map $Nm:p_*(\Gm_S)\to \Gm_{S/\Grp}$.  It induces
a map
$$
Nm: Br(S)[l]=H^2\left((S/\Grp)_{et},p_*(\Gm)\right)[l] \to
H^2\left((S/\Grp)_{et}, \Gm\right)[l]=Br(S/\Grp)[l].
$$
We have $Nm\circ p^*(c)=|\Grp|\cdot c$ and $p^*\circ
Nm(c)=\sum_{\gamma\in \Grp}\gamma^*(c)$. The statement
follows. \endproof

\begin{remark}
This Proposition is suggested by a result of \cite{MR}, where an
analogue of the above Azumaya algebra $\A$ appears. More precisely,
in {\it loc. cit.} one contructs explicitly an Azumaya algebra on
the nilpotent cone in a simple Lie algebra $\g$ over a field of
characteristic $p$.  This algebra is obtained from the enveloping
$U(\g)$ by reduction at the ``singular'' character of the central
subalgebra $Sym({\mathfrak t})^W$; the singular central character
comes from the singular weight $-\rho$ (here ${\mathfrak t}$ is the
maximal torus of $\g$, and $W$ is the Weyl group).
\end{remark}

\begin{theorem}\label{Thm equiv char p}
Assume that the data $\langle V,\omega, \Grp,X \rangle$ over the
positive characteristic residue field $\kk = R'/\m$ satisfy the
assumptions of Theorem~\ref{equivAzu}. Then they also satisfy
Theorem~\ref{evb} (and consequently, Theorem~\ref{main}).
\end{theorem}

\proof{} Let $\A$ be the Azumaya algebra constructed in
Proposition~\ref{untwisting}. We claim that there exist equivalences
\begin{gather}
\Coh(X^\tw)\cong \Coh(X^\tw,\bO \otimes \pi^*(\A^{op}));\label{split 1}\\
D^b\left( \Coh\left (X^\tw,\bO \otimes_{\calo_{V ^\tw }}
\A^{op} \right) \right)\cong D^b\left( \bWbar_{\kk
}^\Grp\otimes_{\calo_{V^\tw
}^\Grp} \A^{op}\fmod\right);\label{D afin untwisted}\\
\bWbar_{\kk }^\Grp\otimes_{\calo_{V^\tw }^\Grp} \A^{op}\fmod\cong
\left( \left(\bWbar_{\kk }\otimes_{\calo_{V^\tw }} \eta^*(
\A^{op})\right) \#\Grp
\right)\fmod;\label{Morita untwist}\\
\Coh^\Grp(V^\tw ) \cong \left( \left( \bWbar_{\kk
}\otimes_{\calo_{V ^\tw}} \eta^*( \A^{op})\right)
\#\Grp\right)\fmod.\label{split 2}
\end{gather}
Here \eqref{split 1}, and \eqref{split 2} follow once we show that
the Azumaya algebra $\bO \otimes \pi^*(\A^{op})$ (respectively, $
\bWbar_{\kk }\otimes_{\calo_{V ^\tw}} \eta^*( \A^{op})$) appearing
in the right-hand side is split (respectively, admits a $\Grp$
equivariant splitting). In each case they split on an open
subvariety by the characterizing property of $\A$; hence they split
on the whole variety by \cite[Corollary IV.2.6]{Milne} (Brauer group
of a regular irreducible scheme injects into the Brauer group of its
generic point); the $\Grp$ equivariant structure on the splitting
bundle for $\bWbar_{\kk }\otimes_{\calo_{V ^\tw}} \eta^*( \A^{op})$
extends from $V_0$ to $V$ since $\codim (V\setminus V_0)\geq 2$.

The equivalence \eqref{Morita untwist} is parallel to the one in
Lemma~\ref{Morita}, and follows from the fact that the idempotent
$e=\frac{1}{|\Grp |}\sum_{\grp\in \Grp}$ generates the unit
two-sided ideal in $(\bWbar_{\kk}\otimes \eta^*(\A^{op})){\#\Grp}$
since it does so in $\bWbar_{\kk}{\#\Grp}$.

Finally, the equivalence \eqref{D afin untwisted} is parallel to the
one in Theorem~\ref{equivAzu}. Namely, by the projection formula we
get
$$
H^0\left(\bO\otimes_{\calo_{V ^\tw}^\Grp} \A^{op}\right) =
\bWbar_{\kk }^\Grp\otimes_{\calo_{V^\tw }^\Grp} \A^{op}.
$$
The projection formula also gives vanishing of the higher cohomology
groups of the sheaf in the left hand side of the last equality. The
algebra in the right-hand side has finite homological dimension by
\eqref{Morita untwist}. Thus \eqref{D afin untwisted} follows from
Proposition~\ref{CYProp}.

Composing the equivalences \eqref{split 1}--\eqref{split 2} we get
the equivalence of the categories of coherent sheaves; now for
$\E$ we can take the image of the structure sheaf
$\calo_V{\#\Grp}\in \Coh^\Grp(V)$ under this equivalence. 

It remains to check this $\E$ is a vector bundle; the equivalence of
categories then shows that it satisfies the conditions of Theorem
\ref{evb}. We use the following more explicit description of $\E$.
Let $\B$ be the $\bO \otimes \pi^*(\A^{op})$ module providing the
Morita equivalence between $\bO$ and $ \pi^*(\A)$; and let $B$ be
the $\bW \otimes \eta^*(\A^{op})$ module providing the Morita
equivalence between $\bW$ and $ \eta^*(\A)$.  It follows from the
definition that the above equivalence sends the coherent sheaf $\B
\in Coh(X^\tw)$ to $B\otimes\kk[\Grp]\in Coh^\Grp(V^\tw)$. The
vector bundle $B$ on $V^\tw$ is trivial by the Suslin-Quillen
Theorem. Thus $\calo\#\Grp$ is a direct summand in $B\otimes
\kk[\Grp]$, so $\E$ is a direct summand in $\B$; in particular, it
is a vector bundle.  \endproof

\begin{remark} One can somewhat simplify
the above explicit description of $\E$ if one is only interested in
the restriction of $\E$ to the formal neighborhood of a fiber of
$\pi$. There the Azumaya algebra $\bO$ splits, and indecomposable
summands of $\E$ are indecomposable summands of a splitting vector
bundle for $\bO$.
\end{remark}

\subsection{Characteristic zero.}

We can now finish the proof of Theorem~\ref{main}. It suffices to
prove Theorem~\ref{evb}.

\begin{lemma}\label{k1k2}
Let $\kk_1\subset \kk_2$ be fields (not necessarily algebraically
closed).  Assume $V$, $\omega$, $\Grp$, $X$ are defined over
$\kk_1$. If the statements of Theorem~\ref{evb} holds
for $V$, $\omega$, $\Grp$, $X$ then it also holds for $V_{\kk_2}$,
$\omega_{\kk_2}$, $\Grp$, $X_{\kk_2}$.  Conversely, if it holds
for $V_{\kk_2}$, $\omega_{\kk_2}$, $\Grp$, $X_{\kk_2}$, then it
also holds for $V_{\kk_1'}$, $\omega_{\kk_1'}$, $\Grp$,
$X_{\kk_1'}$ for some finite extension $\kk_1'$ of $\kk_1$.
\end{lemma}

\proof{} The first statement is obvious. To prove the second one,
consider a finite-type $\kk_1$-algebra $R\subset \kk_2$ such that
both the vector bundle $\E_{\kk_2}$ whose existence is asserted in
Theorem~\ref{evb}, and the isomorphism $\End(\E)\cong \calo_V\#\Grp$
are defined over $R$. Then $\kk_2$ can be taken to be the residue
field of any closed point of $\Spec(R)$. \endproof

\proof[Proof of Theorem~\ref{evb}.] Return to the notation of
Subsection~\ref{ass}; replace $R$ with $R'$ provided by
Lemma~\ref{Morita}. Take an arbitrary closed point $\Spec \kk \in
\Spec R$ such that $R$ is regular in $\Spec \kk$, and the residue
field $\kk = R/\m$ has positive characteristic. Let $\wt{R}$ be the
completion of $R$ with respect to the $\m$-adic topology, and let
$\K'$ be its fraction field. By Theorem~\ref{Thm equiv char p},
Theorem~\ref{evb} holds for $X_\kk$. Moreover, by construction the
vector bundle $\E$ on $X_\kk$ satisfies the assumption of
Proposition~\ref{lifts}. Therefore Theorem~\ref{evb} also holds for
$X_{\K'}$. Applying Lemma~\ref{k1k2}, we see that Theorem~\ref{evb}
holds over a finite extension $F$ of the fraction field of the ring
$R$. Since $F$ is isomorphic to a subfield of $\K$, we see that the
Theorem holds in the original situation.  \endproof

\bigskip

\noindent
{\sc Northwestern University\\
Evanston IL 60208\\
USA\\
\mbox{}\\
and\\
\mbox{}\\
Steklov Math Institute\\
Moscow, USSR\\
\mbox{}\\
{\em E-mail addresses\/}: {\tt bezrukav@math.northwestern.edu}\\
\phantom{{\em E-mail addresses\/}: }{\tt kaledin@mccme.ru}}


\begin{thebibliography}{BMR}

\bibitem[BO]{Ogus} P. Berthelot and A. Ogus, {\it Notes on
crystalline cohomology}, Princeton University Press, Princeton,
N.J.; University of Tokyo Press, Tokyo, 1978.
  
\bibitem[BMR]{MR} R. Bezrukavnikov, I. Mirkovi\'c, and D. Rumynin,
{\em Localization of modules for a semisimple Lie algebra in prime
characteristic}, math.RT/0205144.
 
\bibitem[BKR]{BKR} T. Bridgeland, A. King, and M. Reid, {\em The
McKay correspondence as an equivalence of derived categories},
J. Amer. Math. Soc. {\bf 14} (2001), 535--554.

\bibitem[BK]{BK} A. Bondal and M. Kapranov, {\it Representable
functors, Serre functors, and mutations}, Izv. Ak. Nauk, {\bf 35}
(1990).

%\bibitem[BK]{BK} A.I. Bondal, M.M. Kapranov, {\it Predstavimye
%   funktory, funktory Serra i perestrojki}, Izv. Ak. Nauk, {\bf 53}
%(1989), no. 6, 1183--1205.

\bibitem[BOr]{BO} A. Bondal and D. Orlov, {\it Derived categories of
coherent sheaves}, Proceedings of the International Congress of
Mathematicians, Vol. II (Beijing, 2002), 47--56, Higher Ed. Press,
Beijing, 2002.

\bibitem[G]{Ga} O. Gabber, {\em Some theorems on Azumaya algebras},
The Brauer group (Sem., Les Plans-sur-Bex, 1980), Lecture Notes in
Math. {\bf 844}, Springer, Berlin-New York, 1981.

\bibitem[GV]{SV} G. Gonzalez-Sprinberg and J.-L. Verdier, {\em
Construction g\' eom\' etrique de la correspondance de McKay},
Ann. ENS {\bf 16} (1983), 409--449.
  
\bibitem[H]{Haim} M. Haiman, {\em Hilbert schemes, polygraphs, and
the Macdonald positivity conjecture}, JAMS {\bf 14} (2001),
941--1006.

\bibitem[Ha]{Ha} R. Hartshorne, {\em Residues and duality}, Lecture
notes of a seminar on the work of A. Grothendieck, given at Harvard
1963/64. With an appendix by P. Deligne. Lecture Notes in Math. {\bf
20} Springer-Verlag, Berlin-New York 1966.

\bibitem[IN]{IN} Y. Ito and I. Nakamura, {\em McKay correspondence
and Hilbert schemes}, Proc. Japan Acad. Ser. A Math. Sci. {\bf 72}
(1996), 135--138.
  
\bibitem[Ka1]{Ka} D. Kaledin, {\it Dynkin diagrams and crepant
resolutions of singularities}, math.AG/9903157, to appear in Selecta
Math.

\bibitem[Ka2]{K3} D. Kaledin, {\em Sommese Vanishing for non-compact
manifolds},\\ math.AG/0312271.

\bibitem[K]{Kapranov} M. Kapranov, {\it Noncommutative geometry
based on commutator expansions}, J. Reine Angew. Math. {\bf 505}
(1998), 73--118.

\bibitem[KV]{KV} M. Kapranov and E. Vasserot, {\em Kleinian
singularities, derived categories and Hall algebras}, Math. Ann.
{\bf 316} (2000), 565--576.
  
\bibitem[Ke]{kel} B. Keller, {\em On the cyclic homology of exact
categories}, J. Pure Appl. Algebra {\bf 136} (1999), 1--56.

\bibitem[Ki]{King} A. King, {\it Moduli of representations of
finite-dimensional algebras,} Quart. J. Math. Oxford Ser. (2) {\bf
45} (1994), 515--530.

\bibitem[Kr]{Kr} P.B. Kronheimer, {The construction of ALE spaces as
hyper-K\" ahler quotients}, J. Diff. Geom.  {\bf 29} (1989),
665--683.

\bibitem[Lu]{Lu} G. Lusztig, {\it On quiver varieties}, Adv. Math.
{\bf 136} (1998), 141--182.
  
\bibitem[M]{Milne} J.S. Milne, {\em \'Etale cohomology}, Princeton
Math. Series, {\bf 33}, Princeton U. Press (1980).

\bibitem[N1]{Na98} H. Nakajima, {\em Quiver varieties and 
Kac-Moody algebras}, Duke Math. J. {\bf 91} (1998),
515--560.

\bibitem[N2]{N1} H. Nakajima, {\em Lectures on Hilbert schemes of
points on surfaces}, University Lecture Series, {\bf 18}. American
Mathematical Society, Providence, RI, 1999.

\bibitem[R]{reid} M. Reid, {\em McKay correspondence},
alg-geom/9702016 v3, 1997.

\end{thebibliography}
\end{document}